\documentstyle[12pt,amssymb,amsfonts,epsbox,verbatim]{amsart}
\theoremstyle{plain}
\newtheorem{thm}{Theorem}[section]
\newtheorem{lem}[thm]{Lemma}

\newtheorem{prop}[thm]{Proposition}
\newcommand{\braidrel}[1]{\underset{braid}{\underline{#1}}}
\newcommand{\handlerel}[1]{\underset{handle}{\underline{#1}}}

\theoremstyle{definition}

\theoremstyle{remark}

\begin{document}
\baselineskip=20pt
\title{Action of the mapping class group \\ on a complex of curves 
and a presentation for \\ the mapping class group of a surface} 
\author{Susumu Hirose}
\address{Department of Mathematics, 
Faculty of Science and Engineering, 
Saga University, 
Saga, 840 Japan}
\email{hirose@@ms.saga-u.ac.jp}
\thanks{This research was partially supported by Grant-in-Aid for 
Encouragement of Young Scientists (No. 10740035), Minsitry of Education, 
Science, Sports and Culture, Japan }
\maketitle
\begin{abstract}
Gervais' symmetric presentation for the mapping class group of a surface 
is obtained with using a complex of curves 
in place of Hatcher-Thurston complex. 
\end{abstract}
\section{Introduction}

Let $\Sigma_{g,n}$ be an oriented surface of genus $g$ ($\geq 2$) with 
$n$ ($\geq 0$) boundary components and denote by ${\cal M}_{g,n}$ 
its mapping class group, that is to say, the group of orientation preserving 
diffeomorphisms of $\Sigma_{g,n}$ which are the identity on 
$\partial \Sigma_{g,n}$ modulo isotopy. 
\begin{figure}
\begin{center}
\psbox[height=4cm]{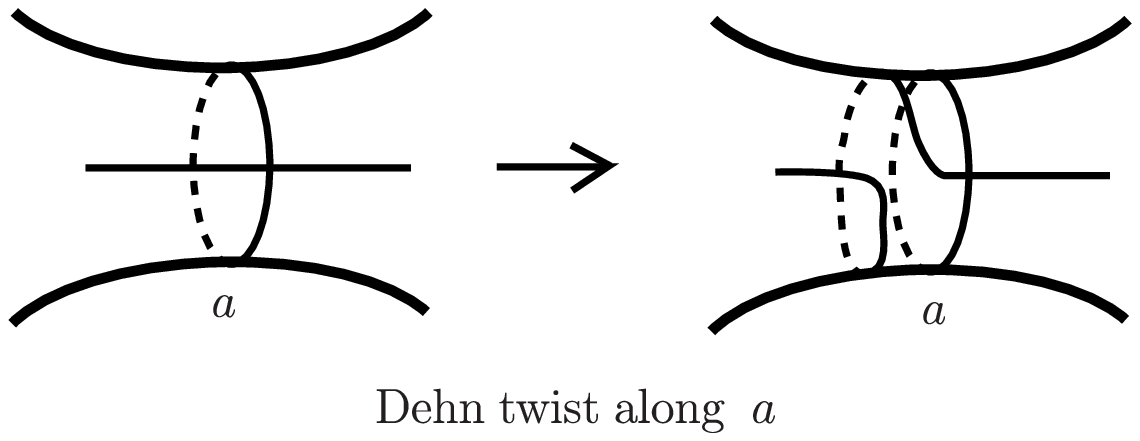}
\caption{}
\label{fig:Dehn-twist}
\end{center}
\end{figure}
For a simple closed curve $a$ in $\Sigma_{g,n}$, we define the Dehn twist along $a$ 
as indicated in Figure \ref{fig:Dehn-twist}. 
We denote the isotopy class of Dehn twist along $a$ by the same letter $a$. 
\par
\begin{figure}
\begin{center}
\psbox[height=7cm]{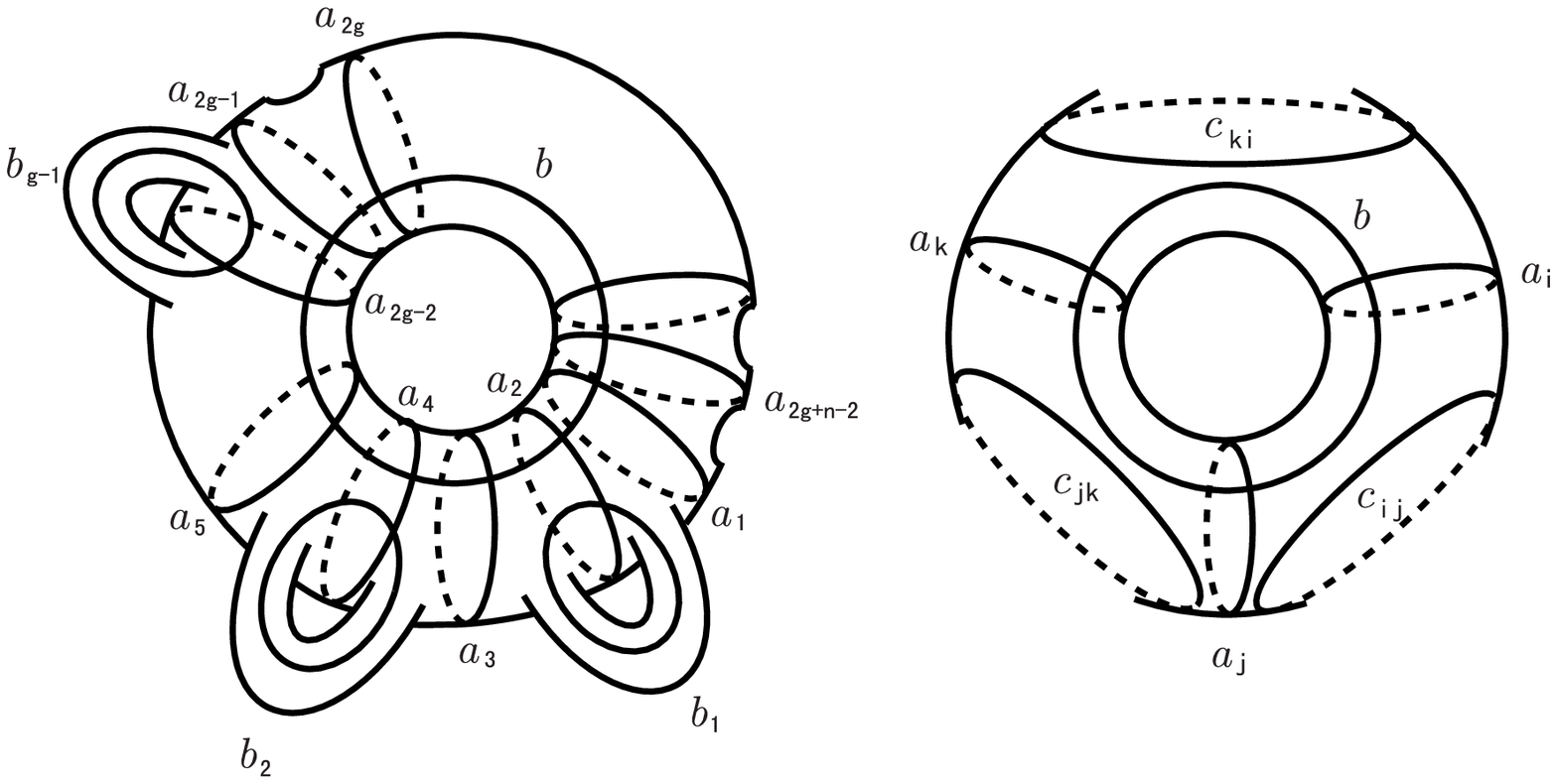}
\caption{}
\label{fig:Gervais}
\end{center}
\end{figure}

It is known that ${\cal M}_{g,n}$ is generated by Dehn twists \cite{Dehn}, 
\cite{Lickorish}. 
McCool \cite{McCool} showed that ${\cal M}_{g,n}$ is finitely presented. 
Hatcher and Thurston \cite{Hatcher-Thurston} defined a simply connected complex 
whose vertices are isotopy classes of "cut systems" and 
introduced a method to give a presentation for ${\cal M}_{g,n}$ 
with using this complex. 
Harer \cite{Harer2} reduced the member of the 2-simplices of this complex, and 
Wajnryb \cite{Wajnryb} gave a simple presentation for ${\cal M}_{g,1}$ and 
${\cal M}_{g,0}$. 
Following from Wajnryb's presentation, Gervais \cite{Gervais} 
gave a symmetric presentation 
for ${\cal M}_{g,n}$. 
We set some notations indicating circles on 
$\Sigma_{g,n}$ as in Figure \ref{fig:Gervais}. 
A triple of integers $(i,j,k)$ $\in$ $\{ 1, \ldots 2g+n-3 \}^3$ will be 
said to be {\it good \/} when: 
\begin{align*}
&i)\ \ (i,j,k) \not\in \{(x,x,x) | x \in \{1.\ldots,2g+n-2\}\}, \\
&ii)\ \ i \leq j \leq k \text{ or } j \leq k \leq i 
\text{ or } k \leq i \leq j. 
\end{align*}
Gervais' symmetric presentation is as follows, 
\begin{thm}\label{thm:Gervais}\cite{Gervais}
If $g \geq 2$, $n \geq 0$, then 
${\cal M}_{g,n}$ is generated by 
$b$, $b_1, \ldots ,b_{g-1}$, $a_1, \ldots, a_{2g+n-2}$, $c_{i,j}$, 
and its defining relations are 
\par\noindent
(A) "HANDLES": $c_{2i, 2i+1} = c_{2i-1, 2i}$ \ for all $i$, $1 \leq i \leq g-1$, 
\par\noindent
(B) "BRAIDS" : for all $x, y$ among the generators, 
$xy = yx$ if the associated curves are disjoint and $xyx = yxy$ if 
the associated curves intersect transversaly in a single point, 
\par\noindent
(C) "STARS"  : 
$c_{ij} c_{jk} c_{ki} = (a_i a_j a_k b)^3$  for all good triples $i,j,k$,   where 
$c_{ll} = 1$. 
\end{thm}
Let $G_{g,n}$ denote the group with presentation given by 
Theorem \ref{thm:Gervais}. 
\par
On the other hand, Harvey \cite{Harvey} 
introduced a complex of curves for $\Sigma_{g,n}$, 
whose vertices are 
isotopy classes of essential (neither homotopic to a point nor 
any boundary component) simple closed curves and 
simplices are the set of vertices 
which are represented by disjoint and non-isotopic curves. 
Harer \cite{Harer} showed the higher connectivity of this complex and, 
with using this complex, the stability of the cohomology group 
of mapping class groups. 
McCullough \cite{McCullough} defined a disk complex of a handle body 
(an oriented 3-dimensional manifold 
obtainted from 3-ball with attaching 1-handles), 
which is defined from a complex of curves with replacing "curves" by 
"meridian disks". 
He showed that the disk complex is contractible. 
The author \cite{Hirose} gave a presentation for the mapping class group of a handle body 
with investigating the action of the mapping class group on this complex. 
The aim of this paper is to give a Gervais' symmetric presentation for 
${\cal M}_{g,n}$ with the same method as above, that is to say, 
with investigating the action of ${\cal M}_{g,n}$ on the complex of curves 
for $\Sigma_{g,n}$. 
We remark here that our method introduced in this paper does not use 
Wajnryb's simple presenation. 
This fact means that we do not need to use Hatcher-Thurston's complex 
to give a presentation for ${\cal M}_{g,n}$. 
\par
Recently, S. Benvenuti (Pisa Univ.) \cite{Benvenuti} showed the similar 
result, independently,  
with using different "complex of curves", which includes separating curves.  
\par
We set notations and conventions used in this paper. 
Composition of elements of ${\cal M}_{g,n}$ will be written from right to left. 
We will denote by $\bar{x}$ 
the inverse of $x$ and $y(x)$ the conjugate $yx\bar{y}$ of $x$ by $y$. 
The notation $\rightleftarrows$ means "commute with", for example, 
for two elements $x$, $y$ of ${\cal M}_{g,n}$, 
$x \rightleftarrows y$ means $xy = yx$. 
We use braid relations and handle relations very often. 
We indicate the place to use a braid relation (resp. handle relation) 
with an underline undersetted by the letter "$braid$" (resp. "$handle$"). 
For example, if $x$, $y$, $z_1$, $z_2$ are loops on $\Sigma_{g,n}$ and if 
$x$ and $y$ intersect transversaly in a single point and 
$z_1$ and $z_2$ are disjoint, then 
$$ \cdots \braidrel{xyx} \cdots \braidrel{z_1 z_2} \cdots 
= \cdots yxy \cdots z_2 z_1 \cdots. $$
\section{A presentation for ${\cal M}_{2,0}$}\label{sec:M_2,0}
Birman and Hilden \cite{Birman-Hilden} showed: 
\begin{thm}\label{thm:Birman-Hilden}\cite{Birman-Hilden} 
${\cal M}_{2,0}$ admits the presentation: \newline
generators: $\tau_1,\ \tau_2,\ \tau_3,\ \tau_4,\ \tau_5$, \newline
defining relations: 
\par
(i) $\tau_i \tau_j = \tau_j \tau_i$, if $|i-j| \geq 2$, $1 \leq i,\ j \leq 5$, 
\par
(ii) $\tau_i \tau_{i+1} \tau_i = \tau_{i+1} \tau_i \tau_{i+1}$ 
\quad $1 \leq i \leq 4$, 
\par
(iii) $(\tau_1 \tau_2 \tau_3 \tau_4 \tau_5)^6 = 1$, 
\par
(iv) $(\tau_1 \tau_2 \tau_3 \tau_4 \tau_5^2 \tau_4 \tau_3 \tau_2 \tau_1)^2 =1$, %
\par
(v) $\tau_1 \tau_2 \tau_3 \tau_4 \tau_5^2 \tau_4 \tau_3 \tau_2 \tau_1 
\rightleftarrows \tau_i$ \quad $1 \leq i \leq 5$. 
\qed
\end{thm}
As we defined previously, 
 $G_{2,0}$ is a group with a following presentation: \newline
generators: $a_1, b, a_2, b_1, c_{1,2}$, \newline
defining relations: \par
(1) $a_1 b a_1 = b a_1 b$, $a_2 b a_2 = b a_2 b$, $a_2 b_1 a_2=b_1 a_2 b_1$, 
$b_1 c_{1,2} b_1 = c_{1,2} b_1 c_{1,2}$, every other pair of generators commutes, 
\par
(2) $(a_1 a_1 a_2 b)^3 = c_{1,2}^2$. 
\par\medskip
Let $\psi_{2,0} : G_{2,0} \rightarrow {\cal M}_{2,0}$ be an epimorphism 
defined by $\psi_{2,0} (a_1) =\tau_1$, $\psi_{2,0} (b) = \tau_2$, 
$\psi_{2,0} (a_2) = \tau_3$, $\psi_{2,0} (b_1) = \tau_4$ and 
$\psi_{2,0} (c_{1,2}) = \tau_5$. 
We want to prove $\psi_{2,0}$ is an isomorphism. 
We shall construct an inverse map 
$\phi_{2,0}$ $: {\cal M}_{2,0} \rightarrow G_{2,0}$. 
For each generators of $G_{2,0}$, we define 
$\phi_{2,0} (\tau_1) = a_1$, $\phi_{2,0} (\tau_2) = b$, 
$\phi_{2,0} (\tau_3) = a_2$, $\phi_{2,0} (\tau_4) = b_1$, 
and $\phi_{2,0} (\tau_5) = c_{1,2}$. 
If the relations (i) - (v) are mapped by $\phi_{2,0}$ onto relations in 
$G_{2,0}$, then $\phi_{2,0}$ extends to homomorphism. 
Then, we can show $\psi_{2,0} \circ \phi_{2,0}$ $= Id_{{\cal M}_{2,0}}$ and 
$\phi_{2,0}$ is an epimorphism, hence, 
$\psi_{2,0}$ is an isomorphism. 
Therefore, in order to prove $\phi_{2,0}$ is an isomorphism, 
it is enough to show that defining relations (i) - (v) are satisfied in 
$G_{2,0}$. 
\par
Relations (i) and (ii) are nothing but the relations (1) for $G_{2,0}$. 
In $G_{2,0}$, the right hand side of relation (v) is 
$a_1 b a_2 b_1 c_{1,2} c_{1,2} b_1 a_2 b a_1$, hence we need to show 
%tmp
\newline
%tmp
$a_1 b a_2 b_1 c_{1,2} c_{1,2} b_1 a_2 b a_1$ 
$\rightleftarrows$ $a_1, b, a_2, b_1, c_{1,2}$. 
For short, we denote $E = a_1 b a_2 b_1 c_{1,2} c_{1,2} b_1 a_2 b a_1$. 
With using the relations (1), 
we can show $E(b) = b$, $E(a_2) = a_2$, $E(b_1) = b_1$, $E(c_{1,2}) = c_{1,2}$. 
In order to show $E(a_1) = a_1$, we have to give 
another presentation for $E$. 
\begin{lem}\label{lem:mirror-star}
$(c_{1,2} c_{1,2} a_2 b_1)^3 = a_1 a_1$. 
\end{lem}
\begin{pf}
We introduce an element 
$D$ $= a_1 b a_2 b_1 c_{1,2} a_1 b a_2 b_1 a_1 b a_2 a_1 b a_1$ 
of ${\cal M}_{2,0}$. 
With using the relations (1), we can show 
$D(a_1) = c_{1,2}$, $D(b) = b_1$, $D(a_2) = a_2$, 
$D(b_1) = b$, and $D(c_{1,2}) = a_1$. 
We take a conjugation of the relation (2) by $D$, then we get 
the equation we need. 
\end{pf}
\begin{lem}\label{lem:another-E} 
$E = a_1 a_1 b a_1 a_1 b \bar{c}_{1,2} \bar{c}_{1,2} \ 
\bar{b}_2 \ \bar{c}_{1,2}\  \bar{c}_{1,2} \ 
\bar{b}_2 $.
\end{lem}
\begin{pf}
By the relations (1), we can show, 
\begin{align*}
c_{1,2} c_{1,2} a_2 b_1 c_{1,2} c_{1,2} a_2 b_1 \braidrel{c_{1,2} c_{1,2} a_2} b_1 &= 
c_{1,2} c_{1,2} a_2 b_1 c_{1,2} c_{1,2} \braidrel{a_2 b_1 a_2} c_{1,2} c_{1,2} b_1 \\
&= c_{1,2} c_{1,2} a_2 b_1 c_{1,2} c_{1,2} b_1 a_2 b_1 c_{1,2} c_{1,2} b_1, 
\end{align*}
We have shown 
$a_1 a_1 = (c_{1,2} c_{1,2} a_2 b_1)^3$, 
in Lemma \ref{lem:mirror-star}, 
therefore, 
$$a_1 a_1 = 
c_{1,2} c_{1,2} a_2 b_1 c_{1,2} c_{1,2} b_1 a_2 b_1 c_{1,2} c_{1,2} b_1.$$ 
From this equation, 
$$
a_2 b_1 c_{1,2} c_{1,2} b_1 a_2 = 
\braidrel{\bar{c}_{1,2} \bar{c}_{1,2} a_1 a_1} \bar{b}_1 \bar{c}_{1,2} 
\bar{c}_{1,2} \bar{b}_1 =
a_1 a_1 \bar{c}_{1,2} \bar{c}_{1,2} \bar{b}_1 \bar{c}_{1,2} \bar{c}_{1,2} \bar{b}_1, 
$$
hence, we can show, 
{\allowdisplaybreaks
\begin{align*}
E &= a_1 b a_2 b_1 c_{1,2} c_{1,2} b_1 a_2 b a_1 
= a_1 b a_1 a_1 \braidrel{\bar{c}_{1,2} \bar{c}_{1,2} \bar{b}_1 
\bar{c}_{1,2} \bar{c}_{1,2} \bar{b}_1 b a_1 } 
\qquad \text{ by the above equation} \\
&=a_1 b a_1 \braidrel{a_1 b a_1} \bar{c}_{1,2} \bar{c}_{1,2} \bar{b}_1 
\bar{c}_{1,2} \bar{c}_{1,2} \bar{b}_1 
= a_1 \braidrel{b a_1 b} a_1 b \bar{c}_{1,2} \bar{c}_{1,2} \bar{b}_1 
\bar{c}_{1,2} \bar{c}_{1,2} \bar{b}_1 \\
&= a_1 a_1 b a_1 a_1 b \bar{c}_{1,2} \bar{c}_{1,2} \bar{b}_1 
\bar{c}_{1,2} \bar{c}_{1,2} \bar{b}_1. 
\end{align*}
}
\end{pf}
We can show $E(a_1) = a_1$ by using the above Lemma and the relations (1).
\par
The relation (iv) is interpreted as $E^2 = 1$ in $G_{2,0}$. 
By Lemma \ref{lem:another-E}, 
{\allowdisplaybreaks
\begin{align*}
E^2 &= a_1 a_1 b a_1 a_1 b \braidrel{\bar{c}_{1,2} \bar{c}_{1,2} \bar{b}_1 
\bar{c}_{1,2} \bar{c}_{1,2} \bar{b}_1 
a_1 a_1 b a_1 a_1 b} \bar{c}_{1,2} \bar{c}_{1,2} \bar{b}_1 
\bar{c}_{1,2} \bar{c}_{1,2} \bar{b}_1 \\
&= a_1 a_1 b a_1 a_1 b a_1 a_1 b a_1 a_1 b 
\bar{c}_{1,2} \bar{c}_{1,2} \bar{b}_1 \bar{c}_{1,2} \bar{c}_{1,2} \bar{b}_1 
\bar{c}_{1,2} \bar{c}_{1,2} \bar{b}_1 \bar{c}_{1,2} \bar{c}_{1,2} \bar{b}_1. 
\end{align*}
}
If we can show $(a_1 a_1 b)^4 = (c_{1,2} c_{1,2} b_1)^4$, then 
$E^2 = (c_{1,2} c_{1,2} b_1)^4 (\bar{c}_{1,2} \bar{c}_{1,2} \bar{b}_2)^4$. 
Since we can show $(c_{1,2} c_{1,2} b_1)^4 = (b_1 c_{1,2} c_{1,2})^4$ 
by the relations (1), we get $E^2 = 1$. 
Therefore it is enough to show: 
\begin{lem}\label{lem:chain-exchange} 
$(a_1 a_1 b)^4 = (c_{1,2} c_{1,2} b_1)^4$ 
\end{lem}
\begin{pf}
We denote $r_1 = a_1 a_1 b a_1 a_1 b$, 
$r_2 = c_{1,2} c_{1,2} b_1 c_{1,2} c_{1,2} b_1$ for short. 
We can show, 
{\allowdisplaybreaks
\begin{align*}
r_1 a_2 r_1 a_2 &= a_1 a_1 b a_1 a_1 b a_2 a_1 \braidrel{a_1 b a_1} a_1 b a_2 
= a_1 a_1 b a_1 a_1 b a_2 \braidrel{a_1 b a_1} b a_1 b a_2 \\
&= a_1 a_1 b a_1 a_1 \braidrel{b a_2 b} a_1 b b a_1 b a_2 
= a_1 a_1 b \braidrel{a_1 a_1 a_2} b \braidrel{a_2 a_1} b b a_1 b a_2 \\
&= a_1 a_1 b a_2 a_1 \braidrel{a_1 b a_1} a_2 b b a_1 b a_2 
= a_1 a_1 b a_2 \braidrel{a_1 b a_1} b a_2 b b a_1 b a_2 \\
&= a_1 a_1 \braidrel{b a_2 b} a_1 b b a_2 b b a_1 b a_2 
= a_1 a_1 a_2 b a_2 a_1 b b a_2 b \braidrel{b a_1 b} a_2 \\
&= a_1 a_1 a_2 b a_2 a_1 b b a_2 \braidrel{b a_1 b} a_1 a_2 
= a_1 a_1 a_2 b a_2 a_1 b b \braidrel{a_2 a_1} b \braidrel{a_1 a_1 a_2} \\
&= a_1 a_1 a_2 b a_2 a_1 b b a_1 \braidrel{a_2 b a_2} a_1 a_1 
= a_1 a_1 a_2 b a_2 a_1 b \braidrel{b a_1 b} a_2 b a_1 a_1 \\
&= a_1 a_1 a_2 b a_2 a_1 \braidrel{b a_1 b} a_1 a_2 b a_1 a_1 
= a_1 a_1 a_2 b \braidrel{a_2 a_1 a_1} b a_1 a_1 a_2 b a_1 a_1 \\
&= (a_1 a_1 a_2 b)^3 a_1 a_1 
\end{align*}
}
and, by the relation (2), 
$(a_1 a_1 a_2 b)^3 a_1 a_1 = c_{1,2} c_{1,2} a_1 a_1$ 
, hence $r_1 a_2 r_1 a_2 = c_{1,2} c_{1,2} a_1 a_1$. 
From the last equation, we can show 
$r_1 ^2 = r_1 \bar{a}_2 \bar{r}_1 c_{1,2} c_{1,2} a_1 a_1 \bar{a}_2$. 
In the same way as above, except using Lemma \ref{lem:mirror-star} 
in place of relation (2), we can show 
$r_2^2 = r_2 \bar{a}_2 \bar{r}_2 c_{1,2} c_{1,2} a_1 a_1 \bar{a}_2$. 
If we can show $r_1(a_2) = r_2(a_2)$, then we can get $r_1^2 = r_2^2$. 
In fact, 
{\allowdisplaybreaks
\begin{align*}
r_2(a_2) &= c_{1,2} \braidrel{c_{1,2} b_1 c_{1,2}} c_{1,2} b_1 (a_2) 
=\braidrel{c_{1,2} b_1 c_{1,2}} \ \braidrel{b_1 c_{1,2} b_1} (a_2) \\
&= b_1 c_{1,2} \braidrel{b_1 c_{1,2} b_1} c_{1,2} (a_2) 
= b_1 c_{1,2} c_{1,2} b_1 \braidrel{c_{1,2} c_{1,2} (a_2)} \\
&= b_1 c_{1,2} c_{1,2} b_1 (a_2)  = b_1 (a_1 a_1 a_2 b)^3 b_1 (a_2) 
\quad \text{ by the relation (2)} \\
&= b_1 a_1 a_1 a_2 b a_1 a_1 a_2 b a_1 a_1 a_2 b \braidrel{b_1 (a_2)} 
= b_1 a_1 a_1 a_2 b a_1 a_1 a_2 b a_1 a_1 \braidrel{a_2 b \bar{a}_2} (b_1) \\
&= b_1 a_1 a_1 a_2 b a_1 a_1 a_2 b a_1 a_1 \bar{b} a_2 \braidrel{b (b_1)} 
= b_1 a_1 a_1 a_2 b \braidrel{a_1 a_1 a_2} b a_1 a_1 \bar{b} a_2 (b_1) \\
&= b_1 a_1 a_1 a_2 b a_2 a_1 \braidrel{a_1 b a_1} a_1 \bar{b} a_2 (b_1)
= b_1 a_1 a_1 a_2 b a_2 a_1 b a_1 \braidrel{b a_1 \bar{b}} a_2 (b_1) \\
&= b_1 a_1 a_1 \braidrel{a_2 b a_2} a_1 b a_1 \bar{a}_1 b \braidrel{a_1 a_2} (b_1) 
= b_1 a_1 a_1 b a_2 b a_1 b b a_2 \braidrel{a_1 (b_1)} \\
&= b_1 a_1 a_1 b a_2 \braidrel{b a_1 b} b a_2 (b_1) 
= b_1 a_1 a_1 b a_2 a_1 \braidrel{b a_1 b} a_2 (b_1) 
= b_1 a_1 a_1 b a_2 a_1 a_1 b \braidrel{a_1 a_2} (b_1) \\
&= b_1 a_1 a_1 b a_2 a_1 a_1 b a_2 \braidrel{a_1 (b_1)} 
= b_1 a_1 a_1 b \braidrel{a_2 a_1 a_1} b a_2 (b_1) 
= b_1 a_1 a_1 b a_1 a_1 \braidrel{a_2 b a_2} (b_1) \\
&= b_1 a_1 a_1 b a_1 a_1 b a_2 \braidrel{b (b_1)} 
= b_1 a_1 a_1 b a_1 a_1 b \braidrel{a_2 (b_1)} 
= b_1 \braidrel{a_1 a_1 b a_1 a_1 b \bar{b}_1} (a_2) \\
&= b_1 \bar{b}_1 a_1 a_1 b a_1 a_1 b (a_2) = a_1 a_1 b a_1 a_1 b (a_2) = r_1(a_2) 
\end{align*} 
}
\end{pf}
The relation (iii) is interpreted as $(a_1 b a_2 b_1 c_{1,2})^6 = 1$. 
If we regard $a_1$, $b$, $a_2$, $b_1$, $c_{1,2}$ as generators of the 6-string 
braid group, namely, $a_1$ as an interchange of 1-st and 2-nd string, 
$b$ as an interchange of 2-nd and 3-rd string and so on, 
then $(a_1 b a_2 b_1 c_{1,2})^6$ is a full twist. 
By investigating a figure of a 6-string full twist, or repeatedly applying 
the relations (1), we can show 
$$
(a_1 b a_2 b_1 c_{1,2})^6 = 
(a_1 b a_2 b_1 c_{1,2})^2 b_1 a_2 b a_1 c_{1,2} b_1 a_2 b (a_2 b_1 c_{1,2})^4. 
$$
By Lemma \ref{lem:mirror-star} , 
{\allowdisplaybreaks
\begin{align*}
a_1 a_1 &= (c_{1,2} c_{1,2} a_2 b_1)^3 \\
&= \braidrel{c_{1,2} c_{1,2} a_2} b_1 c_{1,2} c_{1,2} a_2 b_1 
c_{1,2} c_{1,2} a_2 b_1 
= a_2 c_{1,2} \braidrel{c_{1,2} b_1 c_{1,2}} c_{1,2} a_2 b_1 
c_{1,2} c_{1,2} a_2 b_1 \\
&= a_2 \braidrel{c_{1,2} b_1 c_{1,2}} b_1 c_{1,2} a_2 b_1 
c_{1,2} c_{1,2} a_2 b_1 
= a_2 b_1 c_{1,2} b_1 b_1 \braidrel{c_{1,2} a_2} b_1 
c_{1,2} c_{1,2} a_2 b_1 \\
&= a_2 b_1 c_{1,2} b_1 b_1 a_2 \braidrel{c_{1,2} b_1 c_{1,2}} 
c_{1,2} a_2 b_1 
= a_2 b_1 c_{1,2} b_1 \braidrel{b_1 a_2 b_1} c_{1,2} b_1 
c_{1,2} a_2 b_1 \\
&= a_2 b_1 c_{1,2} \braidrel{b_1 a_2 b_1} a_2 c_{1,2} b_1 
c_{1,2} a_2 b_1 
= a_2 b_1 c_{1,2} a_2 b_1 \braidrel{a_2 a_2 c_{1,2}} b_1 
\braidrel{c_{1,2} a_2} b_1 \\
&= a_2 b_1 c_{1,2} a_2 b_1 c_{1,2} a_2 \braidrel{a_2 b_1 a_2} 
c_{1,2} b_1 
= a_2 b_1 c_{1,2} a_2 b_1 c_{1,2} a_2 b_1 a_2 
\braidrel{b_1 c_{1,2} b_1} \\
&= a_2 b_1 c_{1,2} a_2 b_1 c_{1,2} a_2 b_1 \braidrel{a_2 c_{1,2}} 
b_1 c_{1,2} 
= (a_2 b_1 c_{1,2})^4
\end{align*}
}
therefore, 
{\allowdisplaybreaks
\begin{align*}
(a_1 b a_2 b_1 c_{1,2})^6 &= 
(a_1 b a_2 b_1 c_{1,2})^2 b_1 a_2 b a_1 c_{1,2} b_1 a_2 b a_1 a_1 \\
&= a_1 b a_2 b_1 c_{1,2} a_1 b a_2 \braidrel{b_1 c_{1,2} b_1} 
a_2 b a_1 c_{1,2} b_1 a_2 b a_1 a_1 \\
&= a_1 b a_2 b_1 c_{1,2} \braidrel{a_1 b a_2 c_{1,2}} b_1 
\braidrel{c_{1,2} a_2 b a_1} c_{1,2} b_1 a_2 b a_1 a_1 \\
&= a_1 b a_2 b_1 c_{1,2} c_{1,2} a_1 b \braidrel{a_2 b_1 a_2} b a_1 c_{1,2} 
c_{1,2} b_1 a_2 b a_1 a_1 \\
&= a_1 b a_2 b_1 c_{1,2} c_{1,2} \braidrel{a_1 b b_1} a_2 
\braidrel{b_1 b a_1} c_{1,2} c_{1,2} b_1 a_2 b a_1 a_1 \\
&= a_1 b a_2 b_1 c_{1,2} c_{1,2} b_1 a_1 \braidrel{b a_2 b} a_1 b_1 
c_{1,2} c_{1,2} b_1 a_2 b a_1 a_1 \\
&= a_1 b a_2 b_1 c_{1,2} c_{1,2} b_1 \braidrel{a_1 a_2} b 
\braidrel{a_2 a_1} b_1 c_{1,2} c_{1,2} b_1 a_2 b a_1 a_1 \\
&= a_1 b a_2 b_1 c_{1,2} c_{1,2} b_1 a_2 \braidrel{a_1 b a_1} a_2 
b_1 c_{1,2} c_{1,2} b_1 a_2 b a_1 a_1 \\
&= a_1 b a_2 b_1 c_{1,2} c_{1,2} b_1 a_2 b 
(a_1 b a_2 b_1 c_{1,2} c_{1,2} b_1 a_2 b a_1) a_1  
\end{align*}
}
Previously we have shown that 
$a_1 b a_2 b_1 c_{1,2} c_{1,2} b_1 a_2 b a_1 = E$ 
$\rightleftarrows$ $a_1$, hence, 
$(a_1 b a_2 b_1 c_{1,2})^6$ $=$ 
$(a_1 b a_2 b_1 c_{1,2} c_{1,2} b_1 a_2 b a_1)^2$ $=$ $E^2$ $=$ $1$. 
\section{Elementary relations}\label{sec:elementary}
In this section, we assume $g \geq 3$ or $g=2$, $n \geq 1$ . 
We shall prove some relations in $G_{g,n}$ which are frequently used 
in the following sections. 
The first one is known as "lantern relation", which is proved in 
\cite[Lemma 3]{Gervais}, so we omit the proof here: 
\begin{lem}\label{lem:lantern} For all good triples $(i,j,k)$, 
one has in $G_{g,n}$ the relation, 
$$
(L_{i,j,k}): a_i c_{i,j} c_{j,k} a_k = c_{i,k} a_j X a_j \overline{X} 
= c_{i,k} \overline{X} a_j X a_j, 
$$
where $X = b a_i a_k b$. 
\qed
\end{lem}
The next one is :
\begin{lem}\label{lem:exchange} If $i \not= 2k$, one has in $G_{g,n}$ 
the relation,
\begin{align*}
(X_{i,2k}): 
&(1) \overline{b_k a_{2k} c_{2k-1,2k} b_k}(c_{i,2k}) 
= b a_i a_{2k} b (a_{2k-1}), \\
&(2) b_k a_{2k} c_{2k-1,2k} b_k (c_{i,2k}) 
= \overline{b a_i a_{2k} b }(a_{2k-1}), \\
&(3) \overline{b_k a_{2k} c_{2k-1,2k} b_k}(c_{2k,i}) 
= b a_i a_{2k} b (a_{2k+1}), \\
&(4) b_k a_{2k} c_{2k-1,2k} b_k (c_{2k,i}) 
= \overline{b a_i a_{2k} b }(a_{2k+1}). 
\end{align*}
\end{lem}
\begin{pf}
We will prove (1). Other relations are proved in the same way. 
We denote $X_1 = b_k a_{2k} c_{2k-1,2k} b_k$, 
$X_2 = b a_i a_{2k} b$ for short. Then, 
\begin{align*}
\overline{X_2} \ \overline{X_1} (c_{i,2k}) 
&= \bar{b} \bar{a}_{2k} \bar{a}_i \bar{b} \bar{b}_k \bar{c}_{2k-1,2k} 
\bar{a}_{2k} \braidrel{\bar{b}_k (c_{i,2k})} \\
&= \bar{b} \bar{a}_{2k} \bar{a}_i \bar{b} \bar{b}_k \bar{c}_{2k-1,2k} 
\bar{a}_{2k} c_{i,2k} (b_k) 
\end{align*}
The lantern relation $L_{i,2k-1,2k}$ says
$c_{i,2k} = a_{2k} c_{2k-1,2k} c_{i,2k-1} a_i \bar{a}_{2k-1} 
\overline{X_2} \bar{a}_{2k-1} X_2$. 
Therefore, 
{\allowdisplaybreaks
\begin{align*} 
\overline{X_2} \ \overline{X_1} (c_{i,2k}) 
&= \bar{b} \bar{a}_{2k} \bar{a}_i \bar{b} \bar{b}_k \bar{c}_{2k-1,2k} 
\bar{a}_{2k} a_{2k} c_{2k-1,2k} c_{i,2k-1} a_i 
\bar{a}_{2k-1} \overline{X_2} \bar{a}_{2k-1} X_2 (b_k) \\
&= \bar{b} \bar{a}_{2k} \bar{a}_i \bar{b} \bar{b}_k c_{i,2k-1} a_i 
\bar{a}_{2k-1} \bar{b} \bar{a}_{2k} \bar{a}_i \bar{b} \bar{a}_{2k-1} 
b a_i a_{2k} \braidrel{b (b_k)} \\
&= \bar{b} \bar{a}_{2k} \bar{a}_i \bar{b} \bar{b}_k 
\braidrel{c_{i,2k-1} a_i 
\bar{a}_{2k-1} \bar{b} \bar{a}_{2k} \bar{a}_i \bar{b} \bar{a}_{2k-1} 
b a_i a_{2k}} (b_k) \\
&= \bar{b} \bar{a}_{2k} \bar{a}_i \bar{b} \bar{b}_k 
a_i \bar{a}_{2k-1} \bar{b} \bar{a}_{2k} \bar{a}_i \bar{b} \bar{a}_{2k-1} 
b \braidrel{a_i a_{2k}} \braidrel{c_{i,2k-1}(b_k)} \\
&= \bar{b} \bar{a}_{2k} \bar{a}_i \bar{b} \bar{b}_k 
a_i \bar{a}_{2k-1} \bar{b} \bar{a}_{2k} \bar{a}_i 
\braidrel{\bar{b} \bar{a}_{2k-1} b} 
a_{2k} \braidrel{a_i (b_k)} \\
&= \bar{b} \bar{a}_{2k} \bar{a}_i \bar{b} \bar{b}_k 
\braidrel{a_i \bar{a}_{2k-1}} \bar{b} 
\braidrel{\bar{a}_{2k} \bar{a}_i a_{2k-1}} 
\bar{b} \braidrel{\bar{a}_{2k-1} a_{2k}} (b_k) \\
&= \bar{b} \bar{a}_{2k} \bar{a}_i \bar{b} \bar{b}_k 
\bar{a}_{2k-1} \braidrel{a_i \bar{b} \bar{a}_i} a_{2k-1} 
\braidrel{\bar{a}_{2k} \bar{b} a_{2k}} \braidrel{\bar{a}_{2k-1} (b_k)} \\
&= \bar{b} \bar{a}_{2k} \bar{a}_i \bar{b} \bar{b}_k 
\bar{a}_{2k-1} \bar{b} \bar{a}_i 
\braidrel{b a_{2k-1} b} \bar{a}_{2k} \braidrel{\bar{b} (b_k)} 
= \bar{b} \bar{a}_{2k} \bar{a}_i \bar{b} \bar{b}_k 
\bar{a}_{2k-1} \bar{b} \bar{a}_i 
a_{2k-1} b \braidrel{a_{2k-1} \bar{a}_{2k}}(b_k) \\
&= \bar{b} \bar{a}_{2k} \bar{a}_i \bar{b} \bar{b}_k 
\bar{a}_{2k-1} \bar{b} \bar{a}_i 
a_{2k-1} b \bar{a}_{2k} \braidrel{a_{2k-1}(b_k)} 
= \bar{b} \bar{a}_{2k} \bar{a}_i \bar{b} \bar{b}_k 
\bar{a}_{2k-1} \bar{b} \bar{a}_i 
a_{2k-1} b \braidrel{\bar{a}_{2k}(b_k)} \\
&= \bar{b} \bar{a}_{2k} \bar{a}_i \bar{b} \bar{b}_k 
\braidrel{\bar{a}_{2k-1} \bar{b} \bar{a}_i a_{2k-1} b b_k} (a_{2k}) 
= \bar{b} \bar{a}_{2k} \bar{a}_i \bar{b} \bar{b}_k 
b_k \bar{a}_{2k-1} \bar{b} \braidrel{\bar{a}_i a_{2k-1}} b (a_{2k})\\
&= \bar{b} \bar{a}_{2k} \bar{a}_i \bar{b} 
\braidrel{\bar{a}_{2k-1} \bar{b} a_{2k-1} }
\bar{a}_i b (a_{2k}) 
= \bar{b} \bar{a}_{2k} \bar{a}_i \bar{b} b \bar{a}_{2k-1} 
\braidrel{\bar{b} \bar{a}_i b} (a_{2k})\\
&= \bar{b} \bar{a}_{2k} \bar{a}_i \braidrel{\bar{a}_{2k-1} a_i} \bar{b} 
\braidrel{\bar{a}_i (a_{2k})} 
= \bar{b} \bar{a}_{2k} \bar{a}_i a_i \bar{a}_{2k-1} 
\braidrel{\bar{b} (a_{2k})}
= \bar{b} \bar{a}_{2k} \braidrel{\bar{a}_{2k-1} a_{2k}} (b)\\
&= \bar{b} \bar{a}_{2k} a_{2k} \braidrel{\bar{a}_{2k-1} (b)} 
= \bar{b} b(a_{2k-1}) 
= a_{2k-1}
\end{align*}
}
\noindent
\end{pf}
The third one is known as "chain relation": 
\begin{lem}\label{lem:chain} 
One has in $G_{g,n}$ the relation: 
$$
\{ (c_{2g-2,2g-1})^2 a_{2g-2} b_{g-1} \}^3 = a_{2g-3} a_{2g-1}. 
$$
\end{lem}
\begin{pf} 
We denote 
\begin{align*}
D = c_{2g-2, 2g-1} b_{g-1} &a_{2g-2} b a_{2g-3}
c_{2g-2, 2g-1} b_{g-1} a_{2g-2} b c_{2g-2, 2g-1} b_{g-1} a_{2g-2} \\
& \times c_{2g-2, 2g-1} b_{g-1} c_{2g-2, 2g-1} 
\end{align*}
for short. 
By using braid relations, we can show $D(c_{2g-2,2g-1}) = a_{2g-3}$, 
$D(b_{g-1}) = b$, $D(a_{2g-2}) = a_{2g-2}$, $D(a_{2g-3}) = c_{2g-2,2g-1}$. 
For $D(a_{2g-1})$, 
{\allowdisplaybreaks
\begin{align*} 
D&(a_{2g-1}) \\
&= c_{2g-2, 2g-1} b_{g-1} a_{2g-2} b a_{2g-3}
c_{2g-2, 2g-1} b_{g-1} a_{2g-2} b \\
&\qquad \times
\braidrel{c_{2g-2, 2g-1} b_{g-1} a_{2g-2} 
c_{2g-2, 2g-1} b_{g-1} c_{2g-2, 2g-1} (a_{2g-1})}\\
&= c_{2g-2, 2g-1} b_{g-1} a_{2g-2} b 
\braidrel{a_{2g-3} c_{2g-2, 2g-1} b_{g-1} a_{2g-2}}
 b (a_{2g-1})\\
&= c_{2g-2, 2g-1} b_{g-1} \braidrel{a_{2g-2} b c_{2g-2, 2g-1}} 
b_{g-1} a_{2g-2} a_{2g-3} b (a_{2g-1})\\
&= \braidrel{c_{2g-2, 2g-1} b_{g-1} c_{2g-2, 2g-1}} a_{2g-2} 
\braidrel{b b_{g-1}} a_{2g-2} a_{2g-3} b (a_{2g-1})\\
&= b_{g-1} c_{2g-2, 2g-1} \braidrel{b_{g-1} a_{2g-2} b_{g-1}} 
b a_{2g-2} a_{2g-3} b (a_{2g-1})\\
&= b_{g-1} c_{2g-2, 2g-1} a_{2g-2} b_{g-1} 
\braidrel{a_{2g-2} b a_{2g-2}} a_{2g-3} b (a_{2g-1})\\
&= b_{g-1} c_{2g-2, 2g-1} a_{2g-2} b_{g-1} 
b a_{2g-2} \braidrel{b a_{2g-3} b} (a_{2g-1})\\
&= b_{g-1} c_{2g-2, 2g-1} a_{2g-2} b_{g-1} 
b \braidrel{a_{2g-2} a_{2g-3}} b \braidrel{a_{2g-3} (a_{2g-1})}\\
&= b_{g-1} c_{2g-2, 2g-1} a_{2g-2} b_{g-1} 
b a_{2g-3} a_{2g-2} b (a_{2g-1})\\ 
&= b_{g-1} c_{2g-2, 2g-1} a_{2g-2} b_{g-1} 
\overline{b_{g-1} c_{2g-2, 2g-1} a_{2g-2} b_{g-1}} (c_{2g-2,2g-3}) 
\qquad \text{by } X_{2g-3,2g-2}(3) \\ 
&= c_{2g-2,2g-3}
\end{align*}
}
The star relation $E_{2g-3,2g-3,2g-2}$ of $G_{g,n}$ says: 
\begin{align*}
\{ (a_{2g-3})^2 a_{2g-2} b \}^4 & = \handlerel{c_{2g-3,2g-2}} c_{2g-2,2g-3} \\
& = c_{2g-2, 2g-1} c_{2g-2,2g-3} . 
\end{align*}
We take a conjugation of this equation by $\overline{D}$, then 
we get the equation which we need. 
\end{pf}
\section{A presentation for ${\cal M}_{2,1}$}\label{sec:M_2,1}
\begin{figure}
\begin{center}
\psbox[height=6cm]{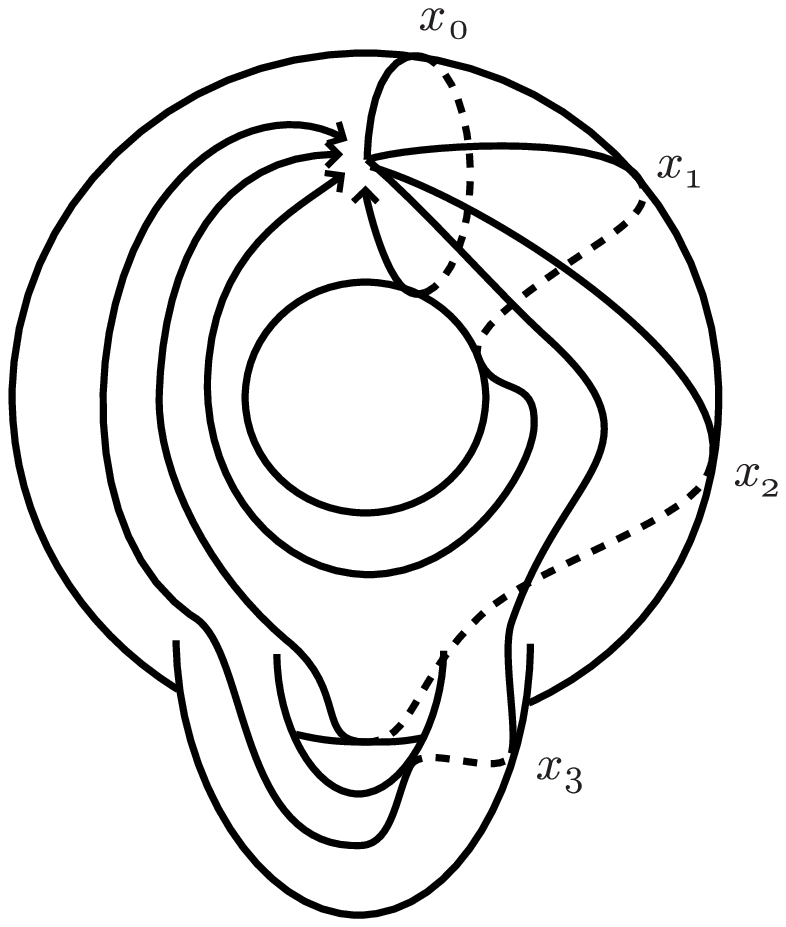}
\caption{}
\label{fig:gen-fund}
\end{center}
\end{figure}
In this section, we give a presentation for ${\cal M}_{2,1}$ and show that 
${\cal M}_{2,1} \cong G_{2,1}$. 
For this purpose, it is enough to show that all the relations for 
${\cal M}_{2,1}$ are satisfied in $G_{2,1}$ by the same reason as section 
\ref{sec:M_2,0}. 
\par
Let $p_1$ be a point on $\Sigma_2$. 
We give a presentation for $\pi_0 (\text{Diff}^+(\Sigma_2, p_1))$ along the 
way of \cite{Birman}. 
Let $\alpha$ be a surjection from $\pi_0 (\text{Diff}^+(\Sigma_2, p_1))$ 
to $\pi_0(\text{Diff}^+(\Sigma_2))$ defined by forgetting a point $p_1$. 
We define a homomorphism $\beta$ from $\pi_1 (\Sigma_2, p_1)$ to 
$\pi_0 (\text{Diff}^+ (\Sigma_2, p_1))$ as follows: 
The homotopy classes of loops indicated in Figure \ref{fig:gen-fund} 
generates $\pi_1(\Sigma_2,p_1)$. 
For a loop $l$ corresponding to one of these generators, 
we take a regular neighborhood $A$ of this loop in $\Sigma_2$. 
Since this $A$ is an annulus, its boundary has two connected components. 
With regarding the orientation for $l$, we denote $A_1$ 
the right hand side of these components, 
and denote $A_2$ the left hand side of them. 
We define 
$\beta($an element of $\pi_1 (\Sigma_2, p_1)$ corresponding to $l)$ $=$ 
$A_1 \bar{A_2}$. 
For short, we denote $x_i = \beta(x_i)$ ($i = 0, 1, 2, 3$). 
For these homomorphisms $\alpha$, $\beta$, there is a short exact sequence: 
\begin{equation}
\tag{S1}
0 \longrightarrow \pi_1 (\Sigma_2, p_1) 
\overset{\beta}{\longrightarrow} \pi_0(\text{Diff}^+(\Sigma_2,p_1)) 
\overset{\alpha}{\longrightarrow} \pi_0(\text{Diff}^+(\Sigma_2)) 
\longrightarrow 0. 
\end{equation}
There is a natural surjection from 
$\pi_0 (\text{Diff}^+(\Sigma_{2,1},rel\ \partial \Sigma_{2,1}))$ to 
%temporary
\newline
$\pi_0 (\text{Diff}^+(\Sigma_{2,1}/\partial \Sigma_{2,1} , 
\partial \Sigma_{2,1}/\partial \Sigma_{2,1}))$ and the latter one is isomorphic to 
$\pi_0(\text{Diff}^+(\Sigma_2, p_1))$, hence there is a surjection $\gamma$ from 
$\pi_0 (\text{Diff}^+(\Sigma_{2,1},rel\ \partial \Sigma_{2,1}))$ 
$\cong$ ${\cal M}_{2,1}$ to 
$\pi_0(\text{Diff}^+(\Sigma_2, p_1))$. 
The kernel of $\gamma$ is an infinite cyclic group ${\Bbb Z}$ generated by 
the Dehn twist along the loop $\partial \Sigma_{2,1}$, 
which we denote $c_{3,1}$. 
Hence, there is a short exact sequence: 
\begin{equation}
\tag{S2}
0 \longrightarrow {\Bbb Z} \longrightarrow {\cal M}_{2,1} 
\overset{\gamma}{\longrightarrow} \pi_0 (\text{Diff}^+(\Sigma_2, p_1)) 
\longrightarrow 0 
\end{equation}
In general, if there is a short exact sequence, 

$$ 0 \longrightarrow L \overset{\phi}{\longrightarrow} 
G \overset{\psi}{\longrightarrow} R \longrightarrow 0, $$
and $L$ and $R$ is finitely presented, then a finite presentation 
for $G$ is given as follows (see, for example, Chapter 10 of \cite{Johnson}). 
Let $l_1, \ldots, l_m$ be the generator of $L$ and, 
$r_1, \ldots, r_n$ be the generator of $R$. 
For each $1 \leq i \leq m$, we denote $\tilde{l_i}$ 
the image of $l_i$ by $\phi$, and 
for each $1 \leq j \leq n$, we fix one of preimages of $r_j$ by $\psi$ and 
denote this $\tilde{r_j}$. 
Then $G$ is generated by $\tilde{l_1}, \ldots, \tilde{l_m}$ and 
$\tilde{r_1}, \ldots, \tilde{r_n}$, and there are following three types of 
relations for $G$. 
\begin{enumerate}
	\item For each $1 \leq i \leq m$, $1 \leq j \leq n$, 
	$\tilde{r_j} \tilde{l_i} \tilde{r_j}^{-1}$ is an element of 
	$\phi(L)$. The equation 
	$$\tilde{r_j} \tilde{l_i} \tilde{r_j}^{-1} = 
	\text{a presentaion of }
	\tilde{r_j} \tilde{l_i} \tilde{r_j}^{-1} 
	\text{ in terms of } \tilde{l_1}, \ldots \tilde{l_m}$$ 
	is a relation for $G$, 
	\item Each relation for $R$ is presented by a word $w(r_1, \ldots, r_n)$. 
	The element 
	%temporaly
	\newline 
	%temporaly
	$w(\tilde{r_1}, \ldots, \tilde{r_n})$ is in the kernel of 
	$\psi$, hence it is an element of $\phi(L)$. 
	The equation 
	$$ w(\tilde{r_1}, \ldots, \tilde{r_n}) = 
	\text{a presentation of } 
	w(\tilde{r_1}, \ldots, \tilde{r_n}) 
	\text{ in terms of } \tilde{l_1}, \ldots \tilde{l_m}$$ 
	is a relation for $G$, 
	\item For each relation for $L$, a equation given from this relation 
	by exchanging $l_i$ with $\tilde{l_i}$ is also a relation for $G$
\end{enumerate}
We apply this method to the above short exact sequences (S1) and (S2). 
For (S1), with remarking that 
$a_1$, $b$, $a_2$, $b_1$, $c_{1,2}$ in $\pi_0(\text{Diff}^+(\Sigma_2, p_1))$ are 
mapped, by $\alpha$, to the elements of $\pi_0(\text{Diff}^+(\Sigma_2))$ denoted 
by the same letters, we can see that 
$\pi_0(\text{Diff}^+(\Sigma_2,p_1))$ is generated by 
$x_0$, $x_1$, $x_2$, $x_3$, $a_1$, $b$, $a_2$, $b_1$, $c_{1,2}$ and 
its defining relations are: 
\begin{equation}\tag{1-$a_1$}
a_1(x_0) = x_0,\ a_1(x_1)=x_1 \bar{x}_0,\ a_1(x_2) = x_2 \bar{x}_0, \ 
a_1(x_3) = x_3 \bar{x}_0, 
\end{equation}
\begin{equation}\tag{1-$b$}
b(x_0) = x_1, \ b(x_1) = x_1 \bar{x}_0 x_1,\ b(x_2) = x_2,\ 
b(x_3) = x_3, 
\end{equation}
\begin{equation}\tag{1-$a_2$}
a_2(x_0) = x_0,\ a_2(x_1) = x_2,\ a_2(x_2)=x_2 \bar{x}_1 x_2,\ 
a_2(x_3) = x_3, 
\end{equation}
\begin{equation}\tag{1-$b_1$}
b_1(x_0) = x_0,\ b_1(x_1) = x_1,\ b_1(x_2)=x_2,\ 
b_1(x_3) = x_3 \bar{x}_2 x_3, 
\end{equation}
\begin{equation}\tag{1-$c_{1,2}$}
c_{1,2} (x_0) = x_0,\ c_{1,2}(x_1) = x_1,\ c_{1,2}(x_2) = x_2,\ 
c_{1,2} (x_3) = x_3 \bar{x}_2 x_1 \bar{x}_0, 
\end{equation}
\begin{equation}\tag{2-1}
\begin{aligned}
&a_1 b a_1 = b a_1 b,\ a_2 b a_2= b a_2 b,\ a_2 b_1 a_2 = b_1 a_2 b_1,\ 
b_1 c_{1,2} b_1 = c_{1,2} b_1 c_{1,2}, \\ 
& \text{other pairs of } \{a_1,\ b,\ a_2,\ b_1,\ c_{1,2} \} 
\text{ commute each other, } 
\end{aligned}
\end{equation}
\begin{equation}\tag{2-2}
(a_1 a_1 a_2 b)^3 \bar{c}_{1,2}^2 \in \beta( \pi_1(\Sigma_2,p_1)), 
\end{equation}
\begin{equation}\tag{3}
x_3 \bar{x}_2 x_1 \bar{x}_0 \bar{x}_3 x_2 \bar{x}_1 x_0 = 1. 
\end{equation}
Among the above relations, 
(1-$a_1$) to (1-$c_{1,2}$) are checked with drawing figures on 
the actions of $a_1$, $b$, $a_2$, $b_1$, $c_{1,2}$ on 
$\pi_1 (\Sigma_2, p_1)$, 
(2-1) and (2-2) come from 
the relation (1) and (2), introduced in section \ref{sec:M_2,0}, 
for ${\cal M}_{2,0}$ $\cong$ $G_{2,0}$, 
(3) is a relation for $\pi_1(\Sigma_2, p_1)$ which is obtained by reading the 
word on the boundary of an octahedron resulting from cutting 
$\Sigma_2$ along $x_0$, $x_1$, $x_2$, $x_3$. 
With using (S2), we can show that ${\cal M}_{2,1}$ is generated by 
$x_0$, $x_1$, $x_2$, $x_3$, $a_1$, $b$, $a_2$, $b_1$, $c_{1,2}$, $c_{3,1}$, 
and the defining relations are 
the relations (1-$a_1$) to (3) up to the powers of $c_{3,1}$. 
On the other hand, we can see 
$x_0 = a_1 \bar{a}_3$, $x_1 = b(x_0)$, $x_2 = a_2(x_1)$, 
$x_3 = b_1(x_2)$, hence, ${\cal M}_{2,1}$ is generated by 
$a_1$, $a_2$, $a_3$, $b$, $b_1$, $c_{1,2}$, $c_{3,1}$. 
We can derive the defining relations for ${\cal M}_{2,1}$ from 
the raltions for $G_{2,1}$ as follows.  
\newline
(1) It is shown, in the proof of Lemma 9 in \cite{Gervais}, 
that all the relations (1-$a_1$) to (1-$c_{1,2}$) up to the powers of 
$c_{3,1}$ are derived from the relations for $G_{2,1}$. 
We remark that 
$$c_{1,2} (x_3) = x_3 \bar{x}_2 x_1 \bar{x}_0 c_{3,1},$$
since this equation shall be used later. 
\newline
(2-1) These relations are nothing but braid relations. 
\newline
(2-2) The lantern relation $L_{2,3,1}$ says 
\begin{align*}
&a_2 c_{2,3} c_{3,1} a_1 = c_{2,1} a_3 X a_3 \bar{X} 
= c_{2,1} \bar{X} a_3 X a_3, \\
&\text{where } X=b a_2 a_1 b ,
\end{align*}
that is to say, 
\begin{align*}
c_{2,1} \bar{c}_{2,3} &= a_2 c_{3,1} a_1 X \bar{a}_3 \bar{X} \bar{a}_3 
\qquad \cdots\cdots (\alpha) \\
a_1 c_{2,3} \bar{a}_3 &= \bar{c}_{3,1} \bar{a}_2 c_{2,1} \bar{X} a_3 X 
\qquad \cdots\cdots (\beta). 
\end{align*}
The star relation $E_{1,1,2}$ says $(a_1 a_1 a_2 b)^3 = c_{1,2} c_{2,1}$, 
therfore, $(a_1 a_1 a_2 b)^3 (\bar{c}_{1,2})^2 = c_{2,1} \bar{c}_{1,2}$. 
For the right hand of the last equation, we can show, 
{\allowdisplaybreaks
\begin{align*}
c_{2,1} \handlerel{\bar{c}_{1,2}} 
&= c_{2,1} \bar{c}_{2,3} 
= \braidrel{a_2 c_{3,1}} a_1 X \bar{a}_3 \bar{X} \bar{a}_3 
\qquad \text{by }(\alpha) \\
&= c_{3,1} a_2 a_1 X \bar{a}_3 \bar{X} \bar{a}_3 
= c_{3,1} a_2 \braidrel{a_1 b a_1} a_2 b \bar{a}_3 \bar{b} \bar{a}_2 
\bar{a}_1 \bar{b} \bar{a}_3 \\
&= c_{3,1} a_2 b a_1 \braidrel{b a_2 b} \bar{a}_3 \bar{b} \bar{a}_2 
\bar{a}_1 \bar{b} \bar{a}_3 
= c_{3,1} a_2 b a_1 a_2 b \braidrel{a_2 \bar{a}_3} \bar{b} \bar{a}_2 
\bar{a}_1 \bar{b} \bar{a}_3 \\
&= c_{3,1} a_2 b a_1 a_2 b \bar{a}_3 \braidrel{a_2 \bar{b} \bar{a}_2} 
\bar{a}_1 \bar{b} \bar{a}_3 
= c_{3,1} a_2 b a_1 a_2 \braidrel{b \bar{a}_3 \bar{b}} \bar{a}_2 
\braidrel{b \bar{a}_1 \bar{b}} \bar{a}_3 \\
&= c_{3,1} a_2 b a_1 \braidrel{a_2 \bar{a}_3} \bar{b} 
\braidrel{a_3 \bar{a}_2 \bar{a}_1} \bar{b} a_1 \bar{a}_3 
= c_{3,1} a_2 b a_1 \bar{a}_3 \braidrel{a_2 \bar{b} \bar{a}_2} \bar{a}_1 
a_3 b a_1 \bar{a}_3 \\
&= c_{3,1} a_2 b a_1 \bar{a}_3 \bar{b} \bar{a}_2 b \bar{a}_1 a_3 b a_1 
\bar{a}_3  
= c_{3,1} x_2 \bar{x}_1 x_0. 
\end{align*}
}
This shows $c_{2,1} \bar{c}_{1,2}$ $\in$ 
$\beta(\pi_1(\Sigma_2, p_1)) \times {\Bbb Z}$. 
Therefore, $(a_1 a_1 a_2 b)^3 (\bar{c}_{1,2})^2$ $\in$ 
$\beta( \pi_1 (\Sigma_2, p_1)) \times {\Bbb Z}$. 
\newline
(3) With using the lantern relation $L_{2,3,1}$ and braid relations, 
we can show, 
{\allowdisplaybreaks
\begin{align*}
x_3 (c_{2,3}) 
&= b_1 a_2 b a_1 \bar{a}_3 \bar{b} \bar{a}_2 \braidrel{\bar{b}_1(c_{2,3})} 
= b_1 a_2 b a_1 \braidrel{\bar{a}_3 \bar{b} \bar{a}_2 c_{2,3}}(b_1) \\
&= b_1 a_2 b a_1 c_{2,3} \bar{a}_3 \bar{b} \bar{a}_2 (b_1) \\
&= b_1 a_2 b \bar{c}_{3,1} \bar{a}_2 c_{2,1} \bar{X} a_3 X 
\bar{b} \bar{a}_2 (b_1) \qquad \text{by }(\beta) \\
&= b_1 a_2 b \bar{c}_{3,1} \bar{a}_2 c_{2,1} \bar{b} 
\braidrel{\bar{a}_1 \bar{a}_2} \bar{b} a_3 b \braidrel{a_2 a_1} b \bar{b} 
\bar{a}_2 (b_1) \\
&= b_1 a_2 b \bar{c}_{3,1} \bar{a}_2 c_{2,1} \bar{b} 
\bar{a}_2 \bar{a}_1 \bar{b} a_3 b a_1 a_2 \bar{a}_2 (b_1) \\
&= b_1 a_2 b \bar{c}_{3,1} \bar{a}_2 c_{2,1} \bar{b} \bar{a}_2 
\braidrel{\bar{a}_1 \bar{b} a_3 b a_1 (b_1)} \\
&= b_1 a_2 b \braidrel{\bar{c}_{3,1} \bar{a}_2 c_{2,1} \bar{b} \bar{a}_2} 
(b_1) 
= b_1 a_2 b \bar{a}_2 c_{2,1} \bar{b} \bar{a}_2 
\braidrel{\bar{c}_{3,1}(b_1)} \\
&= b_1 \braidrel{a_2 b \bar{a}_2 c_{2,1}} \bar{b} \bar{a}_2 (b_1) 
= b_1 c_{2,1} \braidrel{a_2 b \bar{a}_2} \bar{b} \bar{a}_2 (b_1) \\
&= b_1 c_{2,1} \bar{b} a_2 b \bar{b} \bar{a}_2 (b_1) 
= b_1 c_{2,1} \braidrel{\bar{b} (b_1)}
= b_1 \braidrel{c_{2,1}(b_1)} = b_1 \bar{b}_1 (c_{2,1}) 
= c_{2,1}. 
\end{align*}
}
Hence we get: 
{\allowdisplaybreaks
\begin{align*}
c_{2,3}(\bar{x}_3) 
&= c_{2,3} \bar{x}_3 \bar{c}_{2,3} 
= \bar{x}_3 x_3 c_{2,3} \bar{x}_3 \bar{c}_{2,3} \\
&= \bar{x}_3 c_{2,1} \bar{c}_{2,3} \qquad \text{from the above equation } 
x_3(c_{2,3}) = c_{2,1} \\
&= \bar{x}_3 a_2 c_{3,1} a_1 X \bar{a}_3 \bar{X} \bar{a}_3 
\qquad \text{by }(\alpha) \\
&= \bar{x}_3 a_2 \braidrel{c_{3,1} a_1 b a_2 a_1 b \bar{a}_3 \bar{b} \bar{a}_1 
\bar{a}_2 \bar{b} \bar{a}_3} 
= \bar{x}_3 a_2 a_1 b \braidrel{a_2 a_1} b \bar{a}_3 \bar{b} \bar{a}_1 
\bar{a}_2 \bar{b} \bar{a}_3 c_{3,1} \\
&= \bar{x}_3 a_2 \braidrel{a_1 b a_1} a_2 \braidrel{b \bar{a}_3 \bar{b}}
 \bar{a}_1 \bar{a}_2 \bar{b} \bar{a}_3 c_{3,1} 
= \bar{x}_3 a_2 b a_1 b \braidrel{a_2 \bar{a}_3} \bar{b} 
\braidrel{a_3 \bar{a}_1 \bar{a}_2} \bar{b} \bar{a}_3 c_{3,1} \\
&= \bar{x}_3 a_2 b a_1 b \bar{a}_3 \braidrel{a_2 \bar{b} \bar{a}_2} 
a_3 \bar{a}_1 \bar{b} \bar{a}_3 c_{3,1} 
= \bar{x}_3 a_2 b a_1 \braidrel{b \bar{a}_3 \bar{b}} \bar{a}_2 b a_3 
\bar{a}_1 \bar{b} \bar{a}_3 c_{3,1} \\
&= \bar{x}_3 a_2 b a_1 \bar{a}_3 \bar{b} \braidrel{a_3 \bar{a}_2} 
b a_3 \bar{a}_1 \bar{b} \bar{a}_3 c_{3,1} 
= \bar{x}_3 a_2 b a_1 \bar{a}_3 \bar{b} \bar{a}_2 \braidrel{a_3 b a_3} 
\bar{a}_1 \bar{b} \bar{a}_3 c_{3,1} \\
&= \bar{x}_3 a_2 b a_1 \bar{a}_3 \bar{b} \bar{a}_2 b a_3 
\braidrel{b \bar{a}_1 \bar{b}} \bar{a}_3 c_{3,1} 
= \bar{x}_3 a_2 b a_1 \bar{a}_3 \bar{b} \bar{a}_2 b a_3 \bar{a}_1 \bar{b} 
a_1 \bar{a}_3 c_{3,1} \\
&= \bar{x}_3 x_2 \bar{x}_1 x_0 c_{3,1}. 
\end{align*}
}
Previously, we remarked that 
$c_{1,2} (x_3) = x_3 \bar{x}_2 x_1 \bar{x}_0 c_{3,1}$, hence, 
{\allowdisplaybreaks
\begin{align*}
x_3 \bar{x}_2 x_1 \bar{x}_0 \bar{x}_3 x_2 \bar{x}_1 x_0 
&= c_{1,2} x_3 \bar{c}_{1,2} 
\braidrel{ \bar{c}_{3,1} c_{2,3} \bar{x}_3 \bar{c}_{2,3} \bar{c}_{3,1}} \\
&= c_{1,2} x_3 \bar{c}_{1,2} \handlerel{c_{2,3}} \bar{x}_3 
\handlerel{\bar{c}_{2,3}} (\bar{c}_{3,1})^2 \\
&= c_{1,2} x_3 \bar{c}_{1,2} c_{1,2} \bar{x}_3 \bar{c}_{1,2} 
(\bar{c}_{3,1})^2 = (\bar{c}_{3,1})^2 .
\end{align*}
}
This shows that, modulo powers of $c_{3,1}$, 
$x_3 \bar{x}_2 x_1 \bar{x}_0 \bar{x}_3 x_2 \bar{x}_1 x_0 = 1$ is 
derived from relations for $G_{2,1}$. 
\par
From the above results, we conclude: 
\begin{prop}\label{prop:M_2,1} 
${\cal M}_{2,1}$ $\cong$ $G_{2,1}$. \qed
\end{prop}
\section{Action of ${\cal M}_{g,n}$ on $X(\Sigma_{g,n})$ and a presentation for 
${\cal M}_{g,n}$}\label{sec:action}

In this section,  we assume $g \geq 3$, and $ n \geq 1$.  
We call a simple closed curve on $\Sigma_{g,n}$ {\it non-separating\/}, 
if its complement is connected. 
Define a simplicial complex $X(\Sigma_{g,n})$ of dimension $g-1$, 
whose vertices (0-simplices) are the isotopy classes of 
non-separating simple closed curves on $\Sigma_{g,n}$, and whose simplices are 
determined by the rule that a collection of $k+1$ distinct vertices 
spans an $k$-simplex if and only if it admits a collection of representative 
which are pairwise disjoint and the complement of their disjoint union is 
connected. 
This complex $X(\Sigma_{g,n})$ is defined by Harer \cite{Harer}. 
In the same paper, he showed the following Theorem: 
\begin{thm}\label{thm:Harer}\cite[Theorem 1.1]{Harer}
$X(\Sigma_{g,n})$ is homotopy equivalent to a wedge of $(g-1)$-dimensional 
spheres. 
\qed
\end{thm}
Especially, if $g \geq 3$, $X(\Sigma_{g,n})$ is simply connected. 
\par
For each element $\phi$ of ${\cal M}_{g,n}$ and 
a simplex $([C_0], \ldots, [C_n])$ of $X(\Sigma_{g,n})$, 
%temp
\newline
$([\phi(C_0)], \ldots, [\phi(C_n)])$ is also a simplex of 
$X(\Sigma_{g,n})$ . 
Hence, we can define an action of ${\cal M}_{g,n}$ on $X(\Sigma_{g,n})$ 
by $\phi ([C_0], \ldots, [C_n])$ $= ([\phi(C_0)], \ldots, [\phi(C_n)])$. 
We can see that, each of 
$\{ \text{ 2-simplices of } X(\Sigma_{g,n}) \}/ {\cal M}_{g,n}$, 
$\{ \text{ 1-simplices of } X(\Sigma_{g,n}) \}/ {\cal M}_{g,n}$ and 
$\{ \text{ vertices of } X(\Sigma_{g,n}) \}/ {\cal M}_{g,n}$ consists of 
one element, 
each of which is represented by 
$([C_0], [C_1], [C_2])$, $([C_0],[C_1])$, and $([C_0])$, 
where $C_0=c_{2g-2,2g-1}$, $C_1=a_{2g-2}$, $C_2=a_{2g-4}$. 
If the stabilizer of each vertex is finitely presented, 
and if that of each 1-simplex is finitely generated, we can obtain a 
presentation for ${\cal M}_{g,n}$ as in  the way of 
\cite{Laudenbach}, \cite{Wajnryb}. 
Here, we shall recall this method. 
\par
We fix a vertex $v_0$ of $X(\Sigma_{g,n})$, fix an edge 
(= a 1-simplex with orientation) $e_0$ of $X(\Sigma_{g,n})$ 
which emanates from $v_0$, 
and fix a 2-simplex $f_0$ of $X(\Sigma_{g,n})$ which contains $v_0$. 
Let $C_0, C_1$ and $C_2$ be non-separating simple closed curves 
defined as above,
we set $v_0 = [C_0]$, $e_0 = ([C_0], [C_1])$ and $f_0 = ([C_0], [C_1],[C_2])$. 
We choose an element $t_1$ of ${\cal M}_{g,n}$ which switches 
the vertices of $e_0$. 
In our situation, 
we set $t_1 = b_{g-1} c_{2g-2,2g-1} a_{2g-2} b_{g-1}$. 
By this notation, we see $e_0 = (v_0, t_1(v_0))$. 
We denote by $({\cal M}_{g,n})_{v_0}$ the stabilizer of $v_0$, 
by $({\cal M}_{g,n})_{e_0}$ that of $e_0$, and 
by $<t_1>$ an infinite cyclic group generated by $t_1$. 
The free product $({\cal M}_{g,n})_{v_0} * <t_1>$ with 
the following three types of relation 
defines a presentation for ${\cal M}_{g,n}$. 
(In subsection \ref{subsec:vertex}, we give a set of generators for 
$({\cal M}_{g,n})_{v_0}$. 
In the following statements, "a presentation of $s$ as an elements of 
$({\cal M}_{g,n})_{v_0}$" means a presentation of $s$ as a word 
of elements of this set of generators. )
\par
(Y1) $t_1^2 = $ a presentation of $t_1^2$ as an element of 
$({\cal M}_{g,n})_{v_0}$. 
\par
(Y2) For each generator $s$ of $({\cal M}_{g,n})_{e_0}$, 
$$
\begin{aligned}
t_1  
&(\text{ a presentation of } s 
\text{ as an element of } ({\cal M}_{g,n})_{v_0}) 
\overline{t_1} \\
& = \text{ a presentation of } t_1 s \overline{t_1} 
\text{ as an element of } ({\cal M}_{g,n})_{v_0}. 
\end{aligned}
$$
\par
(Y3) For the loop $\partial f_0$ in $X(\Sigma_{g,n})$, we define an element 
$W_{f_0}$ of $({\cal M}_{g,n})_{v_0} * <t_1>$ in the following manner. 
The loop $\partial f_0$ consists of three vertices $v_0, v_1, v_2$ 
and three edges $e_1, e_2, e_3$ 
such that $e_1 = (v_0, v_1)$, $e_2 = (v_1, v_2)$, $e_3 = (v_2, v_0)$.
There is an element $h_1$ of $({\cal M}_{g,n})_{v_0}$ such that 
$h_1(e_0) = e_1$ i.e. $e_1 = (v_0, h_1 t_1(v_0))$, 
then $\overline{h_1 t_1}(e_2)$ is an edge emanating from $v_0$. 
Hence, there is an element $h_2$ of $({\cal M}_{g,n})_{v_0}$ such that 
$h_2(e_0) = \overline{h_1 t_1}(e_2)$ 
i.e. $e_2 = (h_1 t_1 (v_0), h_1 t_1 h_2 t_1 (v_0))$, 
then $\overline{h_1 t_1 h_2 t_1}(e_3)$ is an edge emanating from $v_0$. 
So, there is an element $h_3$ of $({\cal M}_{g,n})_{v_0}$ such that 
$h_3(e_0) = \overline{h_1 t_1 h_2 t_1} (e_3)$ 
i.e. $e_3 = (h_1 t_1 h_2 t_1 (v_0), h_1 t_1 h_2 t_1 h_3 t_1(v_0))$. 
We define $W_{f_0} = h_1 t_1 h_2 t_1 h_3 t_1$. 
This element $W_{f_0}$ fixes $v_0$, so the following is 
a relation for ${\cal M}_{g,n}$: 
$$
W_{f_0} =  \text{ a presentation of } W_{f_0} 
\text{ as an element of } ({\cal M}_{g,n})_{v_0}. 
$$ 
Under the assumption that ${\cal M}_{g-1,n}$ $\cong$ $G_{g-1,n}$, 
if we can show all the relations for $({\cal M}_{g,n})_{v_0}$ and 
the relations of the above three types (Y1) (Y2) (Y3) are satisfied 
in $G_{g,n}$, then we can show the following theorem by the same reason as 
section \ref{sec:M_2,0}. 
\begin{thm}\label{thm:genus-ind}
If $g \geq 3$, $n \geq 1$ and ${\cal M}_{g-1,n} \cong G_{g-1,n}$, then 
${\cal M}_{g,n}$ $\cong$ $G_{g,n}$. 
\end{thm}
In the last section, we have shown ${\cal M}_{2,1}$ $\cong$ $G_{2,1}$ 
(Prop. \ref{prop:M_2,1}), 
therefore, ${\cal M}_{g,1}$ $\cong$ $G_{g,1}$ for any $g \geq 2$. 
On the other hand, Gervais showed the following theorem 
in \S 3 of \cite{Gervais}: 
\begin{thm}\label{thm:boundary} 
If $g\geq 1$, $n \geq 1$ and ${\cal M}_{g,n}$ $\cong$ $G_{g,n}$, 
then ${\cal M}_{g,n+1}$ $\cong$ $G_{g,n+1}$, 
${\cal M}_{g,n-1}$ $\cong$ $G_{g,n-1}$. 
\qed
\end{thm}
Theorem \ref{thm:Gervais} is proved by Theorem \ref{thm:genus-ind} and 
Theorem \ref{thm:boundary}. 
We remark that Theorem \ref{thm:boundary} was proved without using Wajnryb's 
simple presentation \cite{Wajnryb}. 
In the following subsections, we show all relations for 
$({\cal M}_{2,1})_{v_0}$ (subsection \ref{subsec:vertex}), 
relations of type (Y1) and (Y2) (subsection \ref{subsec:edge}), 
and a relation of type (Y3) (subsection \ref{subsec:face}) are 
satisfied in $G_{g,n}$. 
\subsection{A presentation for $({\cal M}_{g,n})_{v_0}$}\label{subsec:vertex}
We assume that ${\cal M}_{g-1,n}$ $\cong$ $G_{g-1,n}$, and $n \geq 1$. 
Let $\text{Diff}^+ (\Sigma_{g,n})$ denote the group of orientation preserving 
diffeomorphisms on $\Sigma_{g,n}$. 
For subsets $A_1, \ldots, A_m$ and $B$ of 
$\Sigma_{g,n}$, we define 
$\text{Diff}^+(\Sigma_{g,n}, A_1, \ldots, A_m, rel \ B)$ $=$ 
$\{ \phi \in \text{Diff}^+(\Sigma_{g,n})\ |\ \phi(A_1)=A_1, \ldots, 
\phi(A_m) = A_m,\ \phi|_B = id_B \}$. 
In this subsection, we give a presentation for 
$ ({\cal M}_{g,n})_{v_0}$ $=$ 
$\pi_0 (\text{Diff}^+ (\Sigma_{g,n}, C_0, rel \ \partial \Sigma_{g,n}))$. 
Let $\Sigma_{g,n}'$ be a surface obtained from $\Sigma_{g,n}$ 
with cutting along $C_0$, 
and let $E_1, E_2$ be connected components of 
$\partial \Sigma_{g,n}'$ which are 
appeared as a result of cutting. 
Let $\alpha$ be a natural surjection from 
$\pi_0 (\text{Diff}^+(\Sigma_{g,n}', E_1 \cup E_2, rel \ 
\partial \Sigma_{g,n}))$ 
to ${\Bbb Z}_2$ which is a permutation group of $E_1$ and $E_2$, and 
$\beta$ be an inclusion of 
$\pi_0 (\text{Diff}^+ (\Sigma_{g,n}', rel \ \partial \Sigma_{g,n}'))$ into 
$\pi_0 (\text{Diff}^+(\Sigma_{g,n}', E_1 \cup E_2, rel \ 
\partial \Sigma_{g,n}))$. 
Then, there is a short exact sequence: 
\begin{align*}
0 \longrightarrow 
& \pi_0 (\text{Diff}^+ (\Sigma_{g,n}', rel \ \partial \Sigma_{g,n}')) 
\overset{\beta}{\longrightarrow} 
\pi_0 (\text{Diff}^+(\Sigma_{g,n}', E_1 \cup E_2, rel \ \partial \Sigma_{g,n})) \\
& \overset{\alpha}{\longrightarrow} 
{\Bbb Z}_2 \longrightarrow 0
\end{align*}
We can see that 
\begin{align*}
\pi_0 (\text{Diff}^+ (\Sigma_{g,n}, C_0, &rel \ \partial \Sigma_{g,n}))\\
&\cong 
\pi_0 (\text{Diff}^+(\Sigma_{g,n}', E_1 \cup E_2, rel \ \partial \Sigma_{g,n}))
/ c_{2g-2,2g-1} = c_{2g-3,2g-2},
\end{align*}
and 
$$\pi_0 (\text{Diff}^+ (\Sigma_{g,n}', rel \ \partial \Sigma_{g,n}')) \cong 
{\cal M}_{g-1,n+2}.
$$ 
By Theorem \ref{thm:boundary}, 
$\pi_0 (\text{Diff}^+ (\Sigma_{g,n}', rel \ \partial \Sigma_{g,n}'))$ $\cong$ 
$G_{g-1,n+2}$. 
Let $r_{g-1}$ $=$ $\{(c_{2g-3,2g-2})^2 b_{g-1} \}^2$, then 
$r_{g-1}$ $\in$ 
$\pi_0 (\text{Diff}^+ (\Sigma_{g,n}, C_0, rel \ \partial \Sigma_{g,n}))$, 
that is to say, we can regard $r_{g-1}$ as an element of 
$\pi_0 (\text{Diff}^+ (\Sigma_{g,n}', E_1 \cup E_2, 
rel\ \partial \Sigma_{g,n}))$. 
Then $\alpha(r_{g-1})$ generates ${\Bbb Z}_2$. 
From the above observations, we can see: 
\par\noindent
$\pi_0 (\text{Diff}^+ (\Sigma_{g,n}, C_0, rel \ \partial \Sigma_{g,n}))$ is 
isomorphic to $G_{g-1,n+2} * <r_{g-1}>$ with the following relations:
\par\noindent
(A1) $c_{2g-2,2g-1} = c_{2g-3,2g-2}$, 
\par\noindent
(A2) For each generator $t$ of $G_{g-1,n+2}$, 
$$
r_{g-1} t \bar{r}_{g-1} = \ \text{ a presentation of } r_g t \bar{r}_g \ 
\text{ as an element of }\ G_{g-1,n+2}, 
$$
\noindent
(A3) $r_{g-1}^2 = c_{2g-3,2g-1}$. 
\par
We need to show that these relations are derived from relations for $G_{g,n}$. 
\newline
(1) The relation (A1) is nothing but a handle relation. 
\newline
(2) By repeatedly applying star relations, we can show 
$G_{g-1,n+2}$ is generated by ${\cal E}$ $=$ 
$\{ b, a_i\ (1 \leq i \leq 2g+n-2), c_{2j-1,2j}\ (1 \leq j \leq g-2),
c_{k-1,k}\ (2g-2 \leq k \leq 2g+n-2),\ c_{2g+n-2,1} \}$. 
Here, we remark that $({\cal M}_{g,n})_{v_0}$ is generated by 
${\cal E}$ $\cup$ $\{ r_{g-1} \}$. 
By drawing figures, we can show: 
{\allowdisplaybreaks
\begin{align*}
& r_{g-1} (b) = b,\ r_{g-1} (a_i) = a_i \quad \text{if } i \not= 2g-2, \  
1 \leq i \leq 2g+n-2\\
& r_{g-1} (c_{2j-1, 2j}) = c_{2j-1, 2j} \quad \text{if } 1 \leq j \leq g-2 \\
& r_{g-1} (c_{2g-3, 2g-2}) = c_{2g-2, 2g-1},\  
r_{g-1}(c_{2g-2,2g-1}) = c_{2g-3,2g-2} \\
&r_{g-1} (c_{k-1,k}) = c_{k-1,k} \quad \text{if } 2g \leq k \leq 2g+n-2 \\
&r_{g-1} (c_{2g+n-2,1}) = c_{2g+n-2,1} \\
&r_{g-1} a_{2g-2} \bar{r}_{g-1} c_{2g-3,2g-1} a_{2g-2} 
= a_{2g-1} a_{2g-3} (c_{2g-3,2g-2})^2 
\ \cdots\cdots (*) 
\end{align*}
}
The above equations except $(*)$ are derived from braid relation. 
We shall show that the equation $(*)$ is satisfied in $G_{g,n}$.  
{\allowdisplaybreaks
\begin{align*}
r_{g-1}(a_{2g-2}) &= 
c_{2g-3,2g-2} \braidrel{c_{2g-3,2g-2} b_{g-1} c_{2g-3,2g-2}} 
c_{2g-3,2g-2} b_{g-1} (a_{2g-2}) \\
&= \braidrel{c_{2g-3,2g-2} b_{g-1} c_{2g-3,2g-2}} \ 
\braidrel{b_{g-1} c_{2g-3,2g-2} b_{g-1}} (a_{2g-2}) \\
&= b_{g-1} c_{2g-3,2g-2} \braidrel{b_{g-1} c_{2g-3,2g-2} b_{g-1}}
c_{2g-3,2g-2} (a_{2g-2}) \\
&= b_{g-1} c_{2g-3,2g-2} c_{2g-3,2g-2} b_{g-1} 
\braidrel{c_{2g-3,2g-2} c_{2g-3,2g-2} (a_{2g-2})} \\
&= b_{g-1} c_{2g-3,2g-2} c_{2g-3,2g-2} b_{g-1} (a_{2g-2})
\end{align*}
}
By a star relation $E_{2g-3,2g-2,2g-1}$ and a handle relation 
$c_{2g-3,2g-2}$ $=$ $c_{2g-2,2g-1}$, 
$$
c_{2g-3,2g-2} c_{2g-3,2g-2} c_{2g-1,2g-3}
= (a_{2g-3} a_{2g-2} a_{2g-1} b)^3,
$$ 
therefore, 
{\allowdisplaybreaks
\begin{align*}
r_{g-1} &(a_{2g-2}) 
= b_{g-1} (a_{2g-3} a_{2g-2} a_{2g-1} b)^3 
\braidrel{\bar{c}_{2g-1,2g-3} b_{g-1}} (a_{2g-2}) \\
&=b_{g-1}(a_{2g-3} a_{2g-2} a_{2g-1} b)^3 b_{g-1} 
\braidrel{\bar{c}_{2g-1,2g-3}(a_{2g-2})} \\
&=b_{g-1}(a_{2g-3} a_{2g-2} a_{2g-1} b)^3 
\braidrel{b_{g-1} (a_{2g-2})}\\
&=b_{g-1}(a_{2g-3} a_{2g-2} a_{2g-1} b)^2 a_{2g-3} 
\braidrel{a_{2g-2} a_{2g-1}} b \bar{a}_{2g-2}(b_{g-1})\\ 
&=b_{g-1}(a_{2g-3} a_{2g-2} a_{2g-1} b)^2 a_{2g-3} 
a_{2g-1} \braidrel{a_{2g-2} b \bar{a}_{2g-2}}(b_{g-1})\\ 
&=b_{g-1}(a_{2g-3} a_{2g-2} a_{2g-1} b)^2 a_{2g-3} 
a_{2g-1} \bar{b} a_{2g-2} \braidrel{b (b_{g-1})}\\ 
&=b_{g-1}(a_{2g-3} a_{2g-2} a_{2g-1} b)^2 a_{2g-3} 
a_{2g-1} \bar{b} \braidrel{a_{2g-2} (b_{g-1})}\\
&=b_{g-1}(a_{2g-3} a_{2g-2} a_{2g-1} b)^2 
\braidrel{a_{2g-3} a_{2g-1} \bar{b} \bar{b}_{g-1}} (a_{2g-2})\\
&=b_{g-1}(a_{2g-3} \braidrel{a_{2g-2} a_{2g-1}} b)^2 
\bar{b}_{g-1} a_{2g-3} a_{2g-1} \bar{b}(a_{2g-2})\\
&=(\braidrel{b_{g-1} a_{2g-3} a_{2g-1}} a_{2g-2} 
\braidrel{b \bar{b}_{g-1}})^2 a_{2g-3} a_{2g-1} \bar{b}(a_{2g-2})\\
&=(a_{2g-3} a_{2g-1} \braidrel{b_{g-1} a_{2g-2} \bar{b}_{g-1}} b)^2 
a_{2g-3} a_{2g-1} \bar{b}(a_{2g-2})\\
&= a_{2g-3} a_{2g-1} \bar{a}_{2g-2} b_{g-1} a_{2g-2} b 
\braidrel{a_{2g-3} a_{2g-1} \bar{a}_{2g-2}} b_{g-1} a_{2g-2} b 
a_{2g-3} a_{2g-1} \bar{b}(a_{2g-2})\\
&= a_{2g-3} a_{2g-1} \bar{a}_{2g-2} b_{g-1} 
\braidrel{a_{2g-2} b \bar{a}_{2g-2}} a_{2g-3} a_{2g-1} b_{g-1} a_{2g-2} b 
a_{2g-3} a_{2g-1} \bar{b}(a_{2g-2})\\
&= a_{2g-3} a_{2g-1} \bar{a}_{2g-2} \braidrel{b_{g-1} \bar{b}} 
 a_{2g-2} \braidrel{b a_{2g-3} a_{2g-1} b_{g-1}} a_{2g-2} b 
a_{2g-3} a_{2g-1} \bar{b}(a_{2g-2})\\
&= a_{2g-3} a_{2g-1} \bar{a}_{2g-2} \bar{b} b_{g-1}
a_{2g-2} b_{g-1} b a_{2g-3} \braidrel{a_{2g-1} a_{2g-2}} b 
\braidrel{a_{2g-3} a_{2g-1}} \bar{b}(a_{2g-2})\\
&= a_{2g-3} a_{2g-1} \bar{a}_{2g-2} \bar{b} b_{g-1}
a_{2g-2} b_{g-1} b a_{2g-3} a_{2g-2} \braidrel{a_{2g-1} b a_{2g-1}}
 a_{2g-3} \bar{b}(a_{2g-2})\\
&= a_{2g-3} a_{2g-1} \bar{a}_{2g-2} \bar{b} b_{g-1}
a_{2g-2} b_{g-1} b a_{2g-3} a_{2g-2} b a_{2g-1} 
\braidrel{b a_{2g-3} \bar{b}}(a_{2g-2})\\
&= a_{2g-3} a_{2g-1} \bar{a}_{2g-2} \bar{b} 
\braidrel{b_{g-1} a_{2g-2} b_{g-1}} b \braidrel{a_{2g-3} a_{2g-2}} 
b \braidrel{a_{2g-1} \bar{a}_{2g-3}} b \braidrel{a_{2g-3}(a_{2g-2})}\\
&= a_{2g-3} a_{2g-1} \braidrel{\bar{a}_{2g-2} \bar{b} a_{2g-2}} b_{g-1} 
\braidrel{a_{2g-2} b a_{2g-2}} \ \braidrel{a_{2g-3} b \bar{a}_{2g-3}}
a_{2g-1} b(a_{2g-2})\\
&= a_{2g-3} a_{2g-1} b \bar{a}_{2g-2} \braidrel{\bar{b} b_{g-1} b} 
a_{2g-2} b \bar{b} a_{2g-3} \braidrel{b a_{2g-1} b}(a_{2g-2})\\
&= a_{2g-3} a_{2g-1} b \bar{a}_{2g-2} b_{g-1}  
a_{2g-2} a_{2g-3} a_{2g-1} b \braidrel{a_{2g-1} (a_{2g-2})}\\
&= a_{2g-3} a_{2g-1} b \bar{a}_{2g-2} b_{g-1}  
a_{2g-2} a_{2g-3} a_{2g-1} \braidrel{b(a_{2g-2})}\\
&= a_{2g-3} a_{2g-1} b \bar{a}_{2g-2} b_{g-1}  
a_{2g-2} \braidrel{a_{2g-3} a_{2g-1} \bar{a}_{2g-2}}(b)\\
&= a_{2g-3} a_{2g-1} b \bar{a}_{2g-2} b_{g-1}  
a_{2g-2} \bar{a}_{2g-2} a_{2g-3} a_{2g-1} (b)\\
&= a_{2g-3} a_{2g-1} b \bar{a}_{2g-2} 
\braidrel{b_{g-1} a_{2g-3} a_{2g-1}} (b)\\
&= a_{2g-3} a_{2g-1} b \braidrel{\bar{a}_{2g-2} a_{2g-3} a_{2g-1}} \ 
\braidrel{b_{g-1} (b)}\\
&= a_{2g-3} a_{2g-1} b a_{2g-3} a_{2g-1} \braidrel{\bar{a}_{2g-2}(b)} 
= \braidrel{a_{2g-3} a_{2g-1}} b a_{2g-3} a_{2g-1} b (a_{2g-2}) \\
&= a_{2g-1} \braidrel{a_{2g-3} b a_{2g-3}} a_{2g-1} b (a_{2g-2}) 
= a_{2g-1} b a_{2g-3} \braidrel{b a_{2g-1} b} (a_{2g-2}) \\
&= a_{2g-1} b \braidrel{a_{2g-3} a_{2g-1}} b \braidrel{a_{2g-1}(a_{2g-2})} 
= \braidrel{a_{2g-1} b a_{2g-1}} a_{2g-3} b (a_{2g-2}) \\
&= b a_{2g-1} \braidrel{b a_{2g-3} b} (a_{2g-2}) 
= b \braidrel{a_{2g-1} a_{2g-3}} b \braidrel{a_{2g-3} (a_{2g-2})} 
= b a_{2g-3} a_{2g-1} b (a_{2g-2}). 
\end{align*}
}
The  lantern relation $L_{2g-3, 2g-2, 2g-1}$ says, 
$$
c_{2g-3,2g-1} a_{2g-2} b a_{2g-3} a_{2g-1} b a_{2g-2} 
\bar{b} \bar{a}_{2g-1} \bar{a}_{2g-3} \bar{b} 
= a_{2g-3} c_{2g-3,2g-2} c_{2g-2,2g-1} a_{2g-1}. 
$$
then, 
{\allowdisplaybreaks
\begin{align*}
b a_{2g-3} a_{2g-1} b a_{2g-2} \bar{b} \bar{a}_{2g-1} \bar{a}_{2g-1} \bar{b} 
&= \braidrel{ \bar{a}_{2g-2} \bar{c}_{2g-3,2g-1} a_{2g-3} 
c_{2g-3,2g-2} c_{2g-2,2g-1} a_{2g-1} } \\
&= a_{2g-1} a_{2g-3} c_{2g-3,2g-2} \handlerel{c_{2g-2,2g-1}} \bar{a}_{2g-2} 
\bar{c}_{2g-3,2g-1} \\
&= a_{2g-1} a_{2g-3} (c_{2g-3,2g-2})^2 \bar{a}_{2g-2} \bar{c}_{2g-3,2g-1} 
\end{align*}
}
Therefore, 
\begin{align*}
b a_{2g-3} a_{2g-1} b a_{2g-2} \bar{b} \bar{a}_{2g-1} \bar{a}_{2g-3} \bar{b} 
c_{2g-3,2g-1} a_{2g-2} 
= a_{2g-1} a_{2g-3} (c_{2g-3,2g-2})^2. 
\end{align*}
In the above equation, we exchange 
$b a_{2g-3} a_{2g-1} b a_{2g-2} 
\bar{b} \bar{a}_{2g-1} \bar{a}_{2g-3} \bar{b}$ 
$=$ $b a_{2g-3} a_{2g-1} b (a_{2g-2})$ with $r_{g-1}(a_{2g-2})$, 
then we get $(*)$. 
Hence the relation (A2) is satisfied in $G_{g,n}$. 
\newline
(3)
At first, we can see: 
{\allowdisplaybreaks
\begin{align*}
r_{g-1} &a_{2g-2} r_{g-1} \braidrel{a_{2g-2} (\bar{c}_{2g-2,2g-1})^2}\\ 
&=\{(\handlerel{c_{2g-3,2g-2}})^2 b_{g-1}\}^2 a_{2g-2} 
\{(\handlerel{c_{2g-3,2g-2}})^2 b_{g-1}\}^2 
\braidrel{a_{2g-2} (\bar{c}_{2g-2,2g-1})^2} \\
&= \{(c_{2g-2,2g-1})^2 b_{g-1}\}^2 a_{2g-2} 
c_{2g-2,2g-1} \braidrel{c_{2g-2,2g-1} b_{g-1} c_{2g-2,2g-1}}
 c_{2g-2,2g-1} b_{g-1}
(\bar{c}_{2g-2,2g-1})^2 a_{2g-2} \\
&= \{(c_{2g-2,2g-1})^2 b_{g-1}\}^2 a_{2g-2} 
\braidrel{c_{2g-2,2g-1} b_{g-1} c_{2g-2,2g-1}} \ 
\braidrel{b_{g-1} c_{2g-2,2g-1} b_{g-1}}
(\bar{c}_{2g-2,2g-1})^2 a_{2g-2} \\
&= \{(c_{2g-2,2g-1})^2 b_{g-1}\}^2 a_{2g-2} 
b_{g-1} c_{2g-2,2g-1} \braidrel{b_{g-1} c_{2g-2,2g-1} b_{g-1}} 
c_{2g-2,2g-1} (\bar{c}_{2g-2,2g-1})^2 a_{2g-2} \\
&= \{(c_{2g-2,2g-1})^2 b_{g-1}\}^2 a_{2g-2} 
b_{g-1} c_{2g-2,2g-1} c_{2g-2,2g-1} b_{g-1} c_{2g-2,2g-1} 
c_{2g-2,2g-1} (\bar{c}_{2g-2,2g-1})^2 a_{2g-2} \\
&= (c_{2g-2,2g-1})^2 b_{g-1} (c_{2g-2,2g-1})^2 
\braidrel{b_{g-1} a_{2g-2} b_{g-1}} (c_{2g-2,2g-1})^2 b_{g-1} a_{2g-2} \\
&= (c_{2g-2,2g-1})^2 b_{g-1} \braidrel{(c_{2g-2,2g-1})^2 a_{2g-2}} 
b_{g-1} a_{2g-2} (c_{2g-2,2g-1})^2 b_{g-1} a_{2g-2} \\
&= \{ (c_{2g-2,2g-1})^2 b_{g-1} a_{2g-2} \}^3 \\
&= a_{2g-3} a_{2g-1} \qquad \text{ by Lemma \ref{lem:chain}}, 
\end{align*}
}
therefore, 
$$
r_{g-1}^2 
= r_{g-1} \bar{a}_{2g-2} \bar{r}_{g-1} a_{2g-1} a_{2g-3} 
(c_{2g-2,2g-1})^2 \bar{a}_{2g-2}. 
$$
From the above equation and $(*)$, we can see 
$r_{g-1}^2 = c_{2g-3,2g-1}$. 
\subsection{Generators of $({\cal M}_{g,n})_{e_0}$, and 
relations of type (Y1) and (Y2)}\label{subsec:edge}
\begin{figure}
\begin{center}
\psbox[height=4cm]{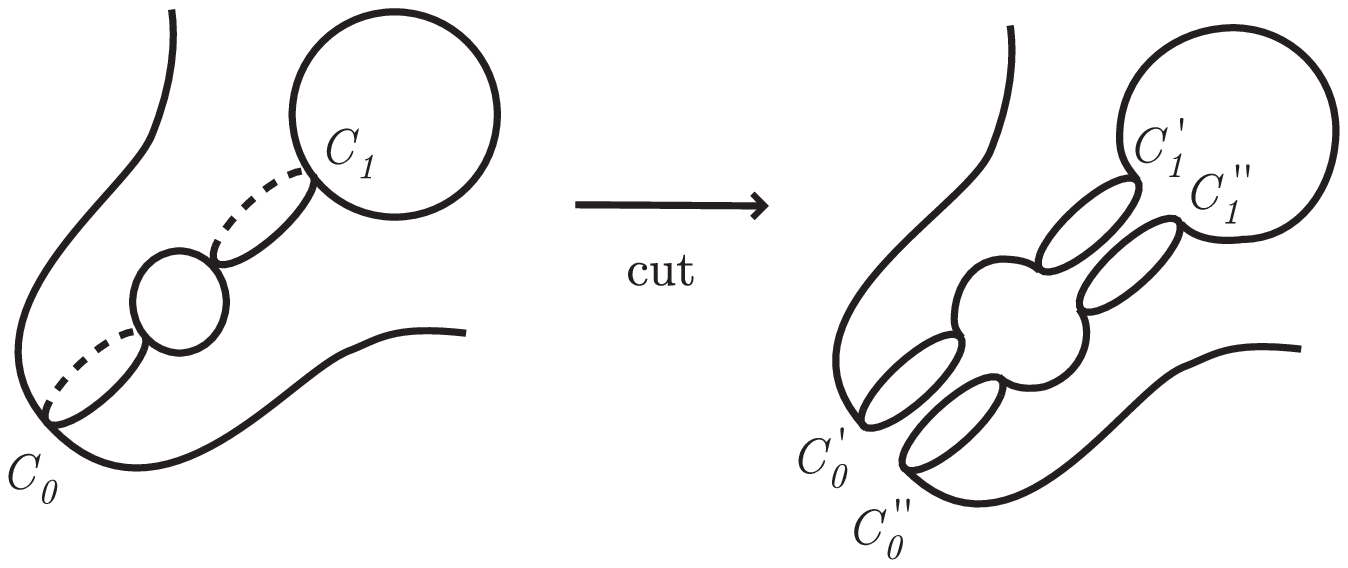}
\caption{}
\label{fig:cutting}
\end{center}
\end{figure}
In this subsection, we give generators of $({\cal M}_{g,n})_{e_0}$ 
$=$ $\pi_0 ( \text{Diff}^+( \Sigma_{g,n}, c_{2g-2,2g-1}, a_{2g-4}, 
rel\ \partial \Sigma_{g,n}))$ and, 
by investigating the action of $t_1$ on these elements, 
we will give relations of type (Y2), and show that 
these relations and a relation of type (Y1) 
are satisfied in $G_{g,n}$. 
\par
At first, we show $t_1^2 \in ({\cal M}_{g,n})_{v_0}$. 
By Lemma \ref{lem:chain} and braid relations, 
{\allowdisplaybreaks
\begin{align*}
&a_{2g-3} a_{2g-1}\\
&= c_{2g-2,2g-1} c_{2g-2,2g-1} a_{2g-2} b_{g-1} 
\braidrel{c_{2g-2,2g-1}c_{2g-2,2g-1}a_{2g-2}} b_{g-1} 
c_{2g-2,2g-1} c_{2g-2, 2g-1} a_{2g-2} b_{g-1} \\
&= c_{2g-2,2g-1} c_{2g-2,2g-1} a_{2g-2} b_{g-1} a_{2g-2} c_{2g-2,2g-1} 
\braidrel{c_{2g-2,2g-1} b_{g-1} c_{2g-2,2g-1}} 
c_{2g-2, 2g-1} a_{2g-2} b_{g-1} \\
&= c_{2g-2,2g-1} c_{2g-2,2g-1} a_{2g-2} b_{g-1} a_{2g-2} 
\braidrel{c_{2g-2,2g-1} b_{g-1} c_{2g-2,2g-1}} 
b_{g-1} c_{2g-2, 2g-1} a_{2g-2} b_{g-1} \\
&= c_{2g-2,2g-1} c_{2g-2,2g-1} a_{2g-2} 
\braidrel{b_{g-1} a_{2g-2} b_{g-1}} 
c_{2g-2,2g-1} b_{g-1} b_{g-1} c_{2g-2, 2g-1} a_{2g-2} b_{g-1} \\
&= c_{2g-2,2g-1} c_{2g-2,2g-1} a_{2g-2} a_{2g-2} b_{g-1} 
\braidrel{a_{2g-2} c_{2g-2,2g-1}} 
b_{g-1} b_{g-1} c_{2g-2, 2g-1} a_{2g-2} b_{g-1} \\
&= (c_{2g-2,2g-1})^2 (a_{2g-2})^2 t_1^2 \qquad
\text{ since } t_1 = b_{g-1} c_{2g-2,2g-1} a_{2g-2} b_{g-1}, 
\end{align*}
}
therefore, 
$t_1^2 = (\bar{a}_{2g-2})^2 (\bar{c}_{2g-2,2g-1})^2 a_{2g-3} a_{2g-1} $ 
$\in ({\cal M}_{g,n})_{v_0}$. 
This shows that the relation of type (Y1) is satisfied in $G_{g,n}$. 
\par
Let $\Sigma_{g,n}''$ be a surface obtained from $\Sigma_{g,n}$ with cutting 
along $C_0 = c_{2g-2,2g-1}$, $C_1 = a_{2g-2}$. 
As in Figure \ref{fig:cutting}, 
let $C_0'$ and $C_0''$ (resp. $C_1'$ and $C_1''$) be connected components of 
$\partial \Sigma_{g,n}''$ which are appeared as a result of cutting along $C_0$ 
(resp. $C_1$). 
We denote the simple closed curve in the interior of $\Sigma_{g,n}''$ which 
is homotopic to $C_0'$ (resp. $C_0''$, $C_1'$, $C_1''$) and 
Dehn twist along this curve by the same letter. 
We can see that $G_{g-2, n+4}$ $\cong$ ${\cal M}_{g-2,n+4}$ is generated by 
$a_i$ $(1 \leq i \leq 2(g-3) + (n+4))$, $b$, 
$b_j$ $(1 \leq j \leq g-3)$, $c_{2k-1,2k}$ $(1 \leq k \leq g-3)$, 
$c_{l,l+1}$ $(2(g-3)+1 \leq l \leq 2(g-3) + (n+3))$, 
and $c_{2(g-3)+(n+4),1}$. 
There is a homomorphism $\gamma$ from ${\cal M}_{g-2,n+4}$ to 
$\pi_0(\text{Diff}^+(\Sigma_{g,n}'', rel \ \Sigma_{g,n}''))$ defined by 
{\allowdisplaybreaks
\begin{align*}
\gamma (a_i) &= c_{i, 2g-4}\quad \text{if } 1 \leq i \leq 2g-5, \\
\gamma (a_{2g-4}) &= c_{2g-4, 2g-2}, \\
\gamma (a_{2g-3}) &= a_{2g-4}, \\
\gamma (a_i) &= c_{i, 2g-4}\quad \text{if } 2g-2 \leq i \leq 2(g-3) + (n+4), \\
\gamma (b) &= b_{g-2}, \\ 
\gamma (b_j) &= b_j\quad \text{if } 1 \leq j \leq g-3, \\   
\gamma (c_{2k-1, 2k}) &= c_{2k-1, 2k}\quad \text{if } 1 \leq k \leq g-3, \\
\gamma (c_{2(g-3)+1, 2(g-3)+2}) &= C_0'', \\
\gamma (c_{2(g-3)+2, 2(g-3)+3}) &= C_1'', \\
\gamma (c_{2(g-3)+3, 2(g-3)+4}) &= C_1', \\
\gamma (c_{2(g-3)+4, 2(g-3)+5}) &= C_0', \\
\gamma (c_{l, l+1}) &= c_{l,l+1} \quad 
\text{if } 2(g-3)+5 \leq l \leq 2(g-3)+(n+3), \\
\gamma (c_{2(g-3)+(n+4), 1}) &= a_{2g-2}. 
\end{align*}
}
This homomorphism is induced by a homeomorphism from $\Sigma_{g-2,n+4}$ 
to $\Sigma_{g,n}''$, hence, $\gamma$ is an isomorphism, and 
this fact means that the set 
\begin{equation*}
{\cal C}_{g,n}'' = 
\left\{ 
\begin{aligned} 
&c_{i,2g-4},\ c_{2g-4,2g-2},\\ 
&b_{g-2},\ b_j,\ c_{2j-1, 2j},\\ 
&C_0',\ C_0'',\ C_1',\ C_1'',\\
& c_{l,l+1},\ c_{2(g-3)+(n+4),1} 
\end{aligned} 
\right. \ 
\left| \ 
\begin{aligned}
&1 \leq i \leq 2g-5,\\ &2g-2 \leq i \leq 2(g-3) + (n+4),\\ 
&1 \leq j \leq g-3, \\ &2(g-3)+5 \leq l \leq 2(g-3)+(n+3) 
\end{aligned}
\right\}
\end{equation*}
generates $\pi_0(\text{Diff}^+(\Sigma_{g,n}'',\ rel \ \Sigma_{g,n}''))$. 
Let ${\Bbb Z}_2 \times {\Bbb Z}_2$ denote the group, whose first factor is 
a permutation group of $C_0'$ and $C_0''$ and second factor is that of 
$C_1'$ and $C_1''$. 
We denote by $\delta$ a natural homomorphism from 
$\pi_0(\text{Diff}^+(\Sigma_{g,n}'', C_0' \cup C_0'', C_1' \cup C_1'', 
rel\ \partial \Sigma_{g,n}))$ to ${\Bbb Z}_2 \times {\Bbb Z}_2$, and 
$\epsilon$ an inclusion of 
$\pi_0(\text{Diff}^+(\Sigma_{g,n}'', rel\ \partial \Sigma_{g,n}''))$ into 
$\pi_0(\text{Diff}^+(\Sigma_{g,n}'', C_0' \cup C_0'', C_1' \cup C_1'', 
rel\ \partial \Sigma_{g,n}))$. 
Then, there is a short exact sequence,
\begin{align*}
0 \longrightarrow \pi_0(&\text{Diff}^+ (\Sigma_{g,n}'', rel \ \Sigma_{g,n}''))\\ 
&\overset{\epsilon}{\longrightarrow} 
\pi_0(\text{Diff}^+(\Sigma_{g,n}'', C_0' \cup C_0'', C_1' \cup C_1'', 
rel\ \partial \Sigma_{g,n}) )
\overset{\delta}{\longrightarrow} {\Bbb Z}_2 \times {\Bbb Z}_2 
\longrightarrow 0
\end{align*}
Let $p=b a_{2g-2} a_{2g-2} b$, $p' = t_1 p \bar{t}_1$, then, 
by drawing some figures, we can check that  
$p$ and $p'$ $\in ({\cal M}_{g,n})_{e_0}$ and 
$p$ (resp. $p'$) reverses the orientation of $C_1$ (resp. $C_0$). 
Hence, $p$ induces a homeomorphism on $\Sigma_{g,n}''$ which 
exchanges $C_0'$ with $C_0''$ (resp. $C_1'$ with $C_1''$). 
On the other hand, there is an isomorphism 
\begin{align*}
\pi_0(\text{Diff}^+( & \Sigma_{g,n}'', C_0' \cup C_0'', C_1' \cup C_1'', 
rel\ \partial \Sigma_{g,n}) )/(C_0' = C_0'', C_1' = C_1'')\\
&\cong
\pi_0 ( \text{Diff}^+(\Sigma_{g,n}, c_{2g-2,2g-1}, a_{2g-2}, 
rel\ \partial \Sigma_{g,n})) ,
\end{align*}
which maps $C_0' = C_0''$ to 
$c_{2g-2,2g-1}$, $C_1' = C_1''$ to $a_{2g-2}$. 
Therefore, we can show that $({\cal M}_{g,n})_{e_0}$ is generated by 
$({\cal C}_{g,n}''- \{ C_0', C_0'', C_1', C_1'' \})$ $\cup $ 
$\{ c_{2g-2,2g-1}, a_{2g-2}, p, p' \}$. 
For each element $s$ of 
${\cal C}_{g,n}''$ $-$ 
$\{c_{2g-2,2g-4}, c_{2g-4,2g-2}, C_0', C_0'', C_1', C_1'' \}$, 
the associated curve of $s$ is disjoint from those of $b_{g-1}$, $a_{2g-2}$, and 
$c_{2g-2,2g-1}$, hence, by braid relations, 
$t_1 s \bar{t}_1$ $=$ $s$ 
$\in ({\cal M}_{g,n})_{v_0}$. 
This fact shows that, for the above element $s$, the relation of type (Y2) is 
satisfied in $G_{g,n}$. 
\par
In subsection \ref{subsec:vertex}, we showed that $({\cal M}_{g,n})_{v_0}$ is generated 
by ${\cal E} \cup \{r_{g-1}\}$, so a presentation of some element as an element of 
$({\cal M}_{g,n})_{v_0}$ means a presentation of this elements as a word of 
${\cal E} \cup \{r_{g-1}\}$. 
Here, we need to present $p$ and $p'$ as a word 
of these elements. 
Since $b, a_{2g-2}$ $\in {\cal E}$, $p$ is presented as an element of 
$({\cal M}_{g,n})_{v_0}$. 
We shall present $p'$ as an element of $({\cal M}_{g,n})_{v_0}$. 
{\allowdisplaybreaks
\begin{align*}
a_{2g-2} b t_1 (b) &= a_{2g-2} b b_{g-1} \braidrel{c_{2g-2,2g-1} a_{2g-2}} \ 
\braidrel{b_{g-1}(b)} \\
&= a_{2g-2} b b_{g-1} a_{2g-2} \braidrel{c_{2g-2,2g-1}(b)} 
= a_{2g-2} b b_{g-1} \braidrel{a_{2g-2} (b)} \\
&= a_{2g-2} \braidrel{b b_{g-1} \bar{b}}(a_{2g-2}) 
= a_{2g-2} \braidrel{b_{g-1} (a_{2g-2})} 
= a_{2g-2} \bar{a}_{2g-2} (b_{g-1}) = b_{g-1}, 
\end{align*}
}
{\allowdisplaybreaks
\begin{align*}
a_{2g-2} b t_1 (a_{2g-2}) &= a_{2g-2} b b_{g-1} c_{2g-2,2g-1} a_{2g-2} 
\braidrel{b_{g-1}(a_{2g-2})} \\
&= a_{2g-2} b b_{g-1} c_{2g-2,2g-1} a_{2g-2} \bar{a}_{2g-2} (b_{g-1}) \\
&= a_{2g-2} b b_{g-1} \braidrel{c_{2g-2,2g-1} (b_{g-1})} 
= a_{2g-2} b b_{g-1} \bar{b}_{g-1} (c_{2g-2,2g-1}) \\
&= \braidrel{a_{2g-2} b (c_{2g-2,2g-1})} = c_{2g-2,2g-1}. 
\end{align*}
}
Here, we remark that this equations shows $t_1(a_{2g-2})$ 
$\in ({\cal M}_{g,n})_{v_0}$. 
From these equations, we can show, 
{\allowdisplaybreaks
\begin{align*}
a_{2g-2} b t_1 p \bar{t}_1 \bar{b} \bar{a}_{2g-2} 
& = a_{2g-2} b t_1 b a_{2g-2} a_{2g-2} b \bar{t}_1 \bar{b} \bar{a}_{2g-2} \\
& = b_{g-1} c_{2g-2,2g-1} c_{2g-2,2g-1} b_{g-1} .
\end{align*}
}
On the other hand, 
\begin{align*}
r_{g-1} &= ((\handlerel{c_{2g-3,2g-2}})^2 b_{g-1})^2 \\
&= ((c_{2g-2,2g-1})^2 b_{g-1})^2 
\end{align*}
Hence, $b_{g-1} c_{2g-2,2g-1} c_{2g-2,2g-1} b_{g-1}$ $=$ 
$(\bar{c}_{2g-1,2g-1})^2 r_{g-1}$. 
From the above equations, we can show 
$p'$ $=$ $t_1 p \bar{t}_1$ $=$ 
$\bar{b} \bar{a}_{2g-1} (\bar{c}_{2g-2,2g-1})^2 r_{g-1} a_{2g-2} b$. 
This gives a presentation of $p'$ as an element of $({\cal M}_{g,n})_{v_0}$. 
For $p$, the relation of type (Y2) is 
\begin{equation*}
t_1 (b a_{2g-2} a_{2g-2} b ) \bar{t}_1 = t_1 p \bar{t}_1 
= \bar{b} \bar{a}_{2g-2} (\bar{c}_{2g-2,2g-1})^2 r_{g-1} a_{2g-2} b 
\end{equation*}
This relation is satisfied in $G_{g,n}$. 
For $p'$, the relation of type (Y2) is, 
\begin{equation*} 
t_1 (\bar{b} \bar{a}_{2g-2} (\bar{c}_{2g-2,2g-1})^2 r_{g-1} a_{2g-2} b) 
\bar{t}_1 \in ({\cal M}_{g,n})_{v_0}. 
\end{equation*}
We shall show that this equation is satisfied in $G_{g,n}$.  
Previously, we have shown $(t_1)^2$, $p$ $\in ({\cal M}_{g,n})_{v_0}$. 
By the definition of $p'$, we can show, 
\begin{equation*}
t_1 (\bar{b} \bar{a}_{2g-2} (\bar{c}_{2g-2,2g-1})^2 r_{g-1} a_{2g-2} b) 
\bar{t}_1 = t_1 (t_1 p \bar{t}_1) \bar{t}_1 = t_1^2 p \bar{t}_1^2 
\in ({\cal M}_{g,n})_{v_0}. 
\end{equation*}
\par
For $c_{2g-2,2g-1}$, $a_{2g-4}$, we can show $t_1$ exchanges 
$c_{2g-2,2g-1}$ and $a_{2g-4}$, 
{\allowdisplaybreaks
\begin{align*}
t_1 (c_{2g-2,2g-1}) 
&= b_{g-1} \braidrel{c_{2g-2,2g-1} a_{2g-1}} \ 
\braidrel{b_{g-1}(c_{2g-2,2g-1})}\\
&= b_{g-1} a_{2g-2} c_{2g-2,2g-1} \bar{c}_{2g-2,2g-1} (b_{g-1}) \\
&= b_{g-1} \braidrel{a_{2g-2}(b_{g-1})} = b_{g-1} \bar{b}_{g-1} (a_{2g-2}) 
= a_{2g-2}, 
\end{align*}
}
{\allowdisplaybreaks
\begin{align*}
t_1 (a_{2g-2}) 
&= b_{g-1} c_{2g-2,2g-1} a_{2g-2} \braidrel{b_{g-1} (a_{2g-2})} \\
&= b_{g-1} c_{2g-2,2g-1} \braidrel{a_{2g-2} \bar{a}_{2g-2}} (b_{g-1}) \\
&= b_{g-1} \braidrel{c_{2g-2,2g-1} (b_{g-1})} 
= b_{g-1} \bar{b}_{g-1} (c_{2g-2,2g-1}) = c_{2g-2,2g-1}. 
\end{align*}
}
This fact shows $t_1 c_{2g-2,2g-1} \bar{t}_1$, $t_1 a_{2g-1} \bar{t}_1$ 
$\in ({\cal M}_{g,n})_{v_0}$. 
\par
For $c_{2g-2,2g-4}$, 
{\allowdisplaybreaks
\begin{align*}
t_1 (c_{2g-2, 2g-4}) 
&= b_{g-1} \braidrel{c_{2g-2,2g-1} a_{2g-2}} b_{g-1} (c_{2g-2,2g-4}) \\
&= b_{g-1} a_{2g-2} \handlerel{c_{2g-2,2g-1}} b_{g-1} (c_{2g-2,2g-4}) \\
&= b_{g-1} a_{2g-2} c_{2g-3,2g-2} b_{g-1} (c_{2g-2,2g-4}) \\
&= \overline{b a_{2g-4} a_{2g-2} b} (a_{2g-1}) 
\qquad (\text{by } X_{2g-4,2g-2} (4)). 
\end{align*}
}
Since $b$, $a_{2g-1}$, $a_{2g-2}$, $a_{2g-4}$ $\in ({\cal M}_{g,n})_{v_0}$, 
this equation shows $t_1 c_{2g-2,2g-4} \bar{t}_1$ 
$\in ({\cal M}_{g,n})_{v_0}$. 
\par
For $c_{2g-4,2g-2}$, we do the same way as above, 
{\allowdisplaybreaks
\begin{align*}
t_1 (c_{2g-4,2g-2}) 
&= b_{g-1} \braidrel{c_{2g-2,2g-1} a_{2g-2}} b_{g-1} (a_{2g-4,2g-2}) \\
&= b_{g-1} a_{2g-2} \handlerel{c_{2g-2,2g-1}} b_{g-1} (c_{2g-4,2g-2}) \\
&= b_{g-1} a_{2g-2} c_{2g-3,2g-2} b_{g-1} (c_{2g-4,2g-2}) \\
&= \overline{b a_{2g-4} a_{2g-2} b} (a_{2g-3}) \qquad 
(\text{by } X_{2g-4,2g-2} (2) ).
\end{align*}
}
Since $b$, $a_{2g-2}$, $a_{2g-3}$, $a_{2g-4}$ 
$\in ({\cal M}_{g,n})_{v_0}$, this equation shows 
$t_1 c_{2g-4,2g-2} \bar{t}_1$ $\in ({\cal M}_{g,n})_{v_0}$. 
\par
Here, we conclude that all the relations of type (Y2) are 
satisfied in $G_{g,n}$. 
\subsection{Relations of type (Y3)}\label{subsec:face}
We define $t_2$ $=$ $b a_{2g-2} a_{2g-4} b$. 
For the notations used to present a relation of type (Y3), 
it is possible to set $h_1 = 1$, $h_2 = t_2$ and $h_3 = t_2$. 
Then, $W_{f_0}$ $=$ $t_1 t_2 t_1 t_2 t_1$. 
By braid relations, we can show $t_1 t_2 t_1$ $=$ $t_2 t_1 t_2$ as follows. 
{\allowdisplaybreaks
\begin{align*}
t_1 t_2(b_{g-1}) 
&= b_{g-1} c_{2g-2,2g-1} a_{2g-2} b_{g-1} b a_{2g-2} 
\braidrel{a_{2g-4} b (b_{g-1})} \\
&= b_{g-1} c_{2g-2,2g-1} a_{2g-2} b_{g-1} b 
\braidrel{a_{2g-2}(b_{g-1})} \\
&= b_{g-1} c_{2g-2,2g-1} a_{2g-2} b_{g-1} 
\braidrel{b \bar{b}_{g-1}} (a_{2g-2}) \\
&= b_{g-1} c_{2g-2,2g-1} a_{2g-2} b_{g-1} \bar{b}_{g-1} b (a_{2g-2}) \\
&= b_{g-1} c_{2g-2,2g-1} a_{2g-2} \braidrel{b (a_{2g-2})} \\
&= b_{g-1} c_{2g-2,2g-1} a_{2g-2} \bar{a}_{2g-2} (b) 
= \braidrel{b_{g-1} c_{2g-2,2g-1} (b)} = b, 
\end{align*}
}
{\allowdisplaybreaks
\begin{align*}
t_1 t_2(c_{2g-2,2g-1}) 
&= b_{g-1} c_{2g-2,2g-1} a_{2g-2} b_{g-1} 
\braidrel{b a_{2g-2} a_{2g-4} b (c_{2g-2,2g-1})} \\
&= b_{g-1} c_{2g-2,2g-1} a_{2g-2} \braidrel{b_{g-1} (c_{2g-2,2g-1})} \\
&= b_{g-1} c_{2g-2,2g-1} \braidrel{a_{2g-2} \bar{c}_{2g-2,2g-1}} (b_{g-1}) \\
&= b_{g-1} c_{2g-2,2g-1} \bar{c}_{2g-2,2g-1} a_{2g-2} (b_{g-1}) \\
&= b_{g-1} \braidrel{a_{2g-2} (b_{g-1})} = b_{g-1} \bar{b}_{g-1} (a_{2g-2}) 
= a_{2g-2}, 
\end{align*}
}
{\allowdisplaybreaks
\begin{align*}
t_1 t_2(a_{2g-2}) 
&= b_{g-1} c_{2g-2,2g-1} a_{2g-2} b_{g-1} b a_{2g-2} a_{2g-4} 
\braidrel{b (a_{2g-2})} \\
&= b_{g-1} c_{2g-2,2g-1} a_{2g-2} b_{g-1} b a_{2g-2} 
\braidrel{a_{2g-4} \bar{a}_{2g-2}} (b) \\
&= b_{g-1} c_{2g-2,2g-1} a_{2g-2} b_{g-1} b a_{2g-2} \bar{a}_{2g-2} 
a_{2g-4} (b) \\
&= b_{g-1} c_{2g-2,2g-1} a_{2g-2} b_{g-1} b \braidrel{a_{2g-4} (b)} \\
&= b_{g-1} c_{2g-2,2g-1} a_{2g-2} b_{g-1} b \bar{b} (a_{2g-4}) \\
&= \braidrel{b_{g-1} c_{2g-2,2g-1} a_{2g-2} b_{g-1} (a_{2g-4})} = a_{2g-4}. 
\end{align*}
}
Therefore, $t_1 t_2 t_1 \bar{t}_2 \bar{t}_1$ $=$ 
$t_1 t_2 (b_{g-1} c_{2g-2,2g-1} a_{2g-2} b_{g-1} ) \bar{t}_2 \bar{t}_1$ 
$=$ $b a_{2g-2} a_{2g-4} b $ $=$ $t_2$, that is 
$t_1 t_2 t_1$ $=$ $t_2 t_1 t_2$. 
So, we get $W_{f_0}$ $=$ $t_1 t_2 t_1 t_2 t_1$ $=$ $t_1^2 t_2 t_1^2$. 
As we have shown in subsection \ref{subsec:edge}, $t_1^2$ 
$\in ({\cal M}_{g,n})_{v_0}$, 
and, since $b$, $a_{2g-2}$, $a_{2g-4}$ $\in ({\cal M}_{g,n})_{v_0}$, 
we can show $t_2 \in ({\cal M}_{g,n})_{v_0}$. 
By using these facts, we conclude that $W_{f_0}$ $\in ({\cal M}_{g,n})_{v_0}$ 
is satisfied in $G_{g,n}$. 

%
%
%%REFERENCES

\end{document}